\documentclass[runningheads]{llncs}

    \makeatletter
    \newcommand\notsotiny{\@setfontsize\notsotiny{5}{6}}
    \makeatother
    
    \let\oldnl\nl
    \newcommand{\nonl}{\renewcommand{\nl}{\let\nl\oldnl}}
    
    \usepackage[hyphens,spaces,obeyspaces]{url}
    
    \newcommand{\cmss}[1]{\textsf{#1}}

    \usepackage{optidef}
    \usepackage[ruled,linesnumbered,noend, vlined]{algorithm2e}
    \SetKwRepeat{Do}{do}{while}
    \usepackage{moreverb}
    \usepackage{cite}
    \usepackage{graphicx}
    \usepackage{booktabs}
    \usepackage{graphics}
    \usepackage{hyperref}
    \usepackage{amsmath}
    \usepackage{amsfonts, amssymb}
    \usepackage{subcaption}
    \usepackage[authoryear]{natbib}
    \usepackage{doi}
    \usepackage{pgfplots}
    \pgfplotsset{compat=1.18}
    \usepackage{pgf-pie}
    \usepackage{tikz}
    \usepackage{sidecap}
    \usepackage{orcidlink} 
    \usepackage[dvipsnames]{xcolor}
    \usetikzlibrary{patterns}
    \usetikzlibrary{patterns.meta}
    \usetikzlibrary{shapes.geometric, arrows}

\pgfdeclarepatternformonly{plus}{\pgfpoint{-4pt}{-4pt}}{\pgfpoint{4pt}{4pt}}{\pgfpoint{6pt}{6pt}}{%
  \pgfsetlinewidth{0.6pt}
  \pgfpathmoveto{\pgfpoint{-2pt}{0pt}}%
  \pgfpathlineto{\pgfpoint{2pt}{0pt}}%
  \pgfpathmoveto{\pgfpoint{0pt}{-2pt}}%
  \pgfpathlineto{\pgfpoint{0pt}{2pt}}%
  \pgfusepath{stroke}%
}

\pgfdeclarepatternformonly{minus}{\pgfpoint{-4pt}{-4pt}}{\pgfpoint{4pt}{4pt}}{\pgfpoint{6pt}{6pt}}{%
  \pgfsetlinewidth{0.6pt}
  \pgfpathmoveto{\pgfpoint{-2pt}{0pt}}%
  \pgfpathlineto{\pgfpoint{2pt}{0pt}}%
  \pgfusepath{stroke}%
}
    
    \tikzstyle{startstop} = [ellipse, draw, text width=4cm, align=center]
    \tikzstyle{process} = [rectangle, draw, text width=4cm, align=center]
    \tikzstyle{decision} = [diamond, draw, text width=3cm, align=center]
    \tikzstyle{arrow} = [thick,->,>=stealth]
    
    \usepackage{mathtools}
    
    \newtagform{brackets}{[}{]}
    \usetagform{brackets}

    \usepackage{wrapfig} 
    \usepackage{xcolor} 
    \usepackage{float} 
    \usepackage[font=small,labelfont=bf,justification=centering]{caption} 
    
    \newcommand\xqed[1]{%
      \leavevmode\unskip\penalty9999 \hbox{}\nobreak\hfill
      \quad\hbox{#1}}

    \DeclareMathAlphabet\mathbfcal{OMS}{cmsy}{b}{n}

    \usepackage{ulem}
    \usepackage[colorinlistoftodos,prependcaption,textsize=tiny]{todonotes}
    
    \usepackage{xr}

    \usepackage[table]{xcolor}
    \usetikzlibrary{shapes, arrows.meta, positioning, calc}

\newcommand{\pp}{\cmss{p}}
\newcommand{\xx}{\cmss{x}}
\newcommand{\rr}{\cmss{r}}
\newcommand{\yy}{\cmss{y}}

\usepackage{dcolumn}    

\newcolumntype{M}{D{|}{|}{-1}}

\begin{document}
\title{Branch and price for nonlinear production-maintenance scheduling in complex machinery}

\author{Jo\~{a}o Dion\'{\i}sio\inst{1,2,3}\thanks{Corresponding author. Email: \texttt{joao.goncalves.dionisio@gmail.com}}\,\orcidlink{0009-0005-5160-0203} \and Ambros Gleixner\inst{3,4}\,\orcidlink{0000-0003-0391-5903} \and Jo\~{a}o Pedro Pedroso\inst{1,2}\,\orcidlink{0000-0003-1298-7191} \and Ksenia Bestuzheva\inst{3}\,\orcidlink{0000-0002-7018-7099}}
\authorrunning{Dionísio \and Gleixner \and Pedroso \and Bestuzheva} 
\institute{Faculdade de Ci\^{e}ncias, Universidade do Porto, rua do Campo Alegre s/n, $4169$–$007$ Porto, Portugal\\
\and CMUP, Centro de Matem\'{a}tica da Universidade do Porto
\and Zuse Institute Berlin, Takustr. $7$, $14195$ Berlin, Germany
\and Hochschule f\"{u}r Technik und Wirtschaft Berlin, Germany
}

\maketitle

\begin{abstract}
This paper proposes a mixed-integer nonlinear programming approach for joint scheduling of long-term maintenance decisions and short-term production for groups of complex machines with multiple interacting components.
We introduce an abstract model where the production and the condition of machines are described by convex functions, allowing the model to be employed for various application areas fitting the scheme.
We develop a branch-and-price algorithm to solve this problem, enhanced with acceleration techniques to find primal solutions and reduce the number of pricing rounds.
An experimental comparison of this approach to solving the compact formulation directly demonstrates the benefit of the decomposition approach, in particular in larger instances.

\end{abstract}

\label{sec:intro}
\section{Introduction}

Industrial maintenance encompasses a wide range of actions and strategies taken to ensure that machinery used in industrial contexts remains functional.
The worldwide cost of industrial maintenance was evaluated at over 54 billion USD in 2024 and is projected to increase to 73 billion USD by 2029~\citep{Industrial_maintenace_services_global_market_report_2024}, driven by multiple factors, notably the increase in global demand.
Efficient, timely maintenance is of particularly high importance in critical infrastructures—such as power plants and water treatment facilities—especially when the underlying equipment is costly and generally old.
These factors motivate the development of new models and methods for maintenance planning.
Since machine degradation and, thus, the need for maintenance actions typically depend on the machines' production output, approaches that integrate production and maintenance planning are of interest.

The problem we tackle can be summarized as follows. 
A group of machines produces a product that must satisfy a predetermined demand, and each machine is composed of multiple interacting components.
The machines' components may interact in the maintenance sense, where the maintenance of one component implies that of another, and also in the degradation sense, where the condition of a component impacts the degradation speed of another.
The production of a machine entails a degradation process, which may trigger machine maintenance, potentially necessary for demand satisfaction at a later time period.
Maintenance actions are applied to the components directly, and the evolution of component condition can depend nonlinearly on both the production and the condition of other components.

Example applications include the minimization of maintenance costs of power transformer fleets, where the ``product'' is electricity, or of water treatment stations, where the ``product'' is clean water.
More abstractly, it can also be understood as production-maintenance scheduling in multiple factories working in parallel, each containing different interacting machines.
In this last example, the factories take on the role of the (potentially different) machines and the machines of the (potentially interacting) components.
Example~\ref{ex:problem_visualization} below illustrates the problem in an abstract setting, and the limitations of disjoint production and maintenance planning.

\begin{example}
Figure~\ref{fig:visual_representation} represents an instance with three time periods and three machines.
Machines $1$ and $2$ each have four components, and machine $3$ has 6 components and a higher potential output, illustrated by its increased complexity.
The arrows represent maintenance implications, i.e., if component $K$ is maintained, then so are components $L$ and $M$.
For simplicity, the components' condition is discrete in this example--green, yellow, orange, red.
To allow for a grayscale interpretation, the component's condition is also number-coded, with green corresponding to 1, yellow to 2, orange to 3, and red to 4.

In this example, we assume that the demand is the same for every period.
The shapes linking machines and demand represent the machines' production.
So, in period $1$, machine $1$ produces roughly half the demand, while machine $2$ produces the least of the three machines.
Each machine's production is visually represented by a different pattern:
machine $1$ with dots, machine $2$ with horizontal lines, and machine $3$ with vertical lines.
The demand is satisfied if the production of all machines covers it.

\input{compact_example.tex}

In the first period, machine 1 produces over half of the total demand.
This increased production leads to significant degradation of machine $1$'s components, especially component A, which must be maintained.
Component C will also be maintained, despite its degradation being only partial, because of the A$\rightarrow$C maintenance implication.
As machines cannot produce while under maintenance, the other machines must compensate in period 2.
This compensation is very damaging to machine 2, necessitating maintenance on all of its components for it to become functional again.
This forces machines 1 and 3 to operate at their maximum capacity in the third period.
Their maximum capacity is reduced due to the damage to their components---a degraded machine is not able to produce as much as one in optimal condition.
Even when operating at current maximum capacity, the machines cannot satisfy all the demand, and thus this solution is not feasible.

Two suboptimal decisions in this solution are the over-production by machine 1 in the first period and the non-replacement of components B and D during the second period.
Even though these components did not strictly require maintenance as they were still operational, their damaged state did not allow machine 1 to operate at maximum capacity in the third period, which would have been sufficient to satisfy the demand.
These limitations of separating production and maintenance decisions motivate our work on integrated approaches.

\xqed{$\triangle$}
\label{ex:problem_visualization}
\end{example}

\subsection{Related Work}

\subsubsection{Maintenance in practice} Industrial maintenance tends to rely on two heuristics for maintenance scheduling: Time-Based Maintenance (TbM), where maintenance is planned considering the time since the last maintenance action, and Condition-Based Maintenance (CbM), where the asset is regularly tested and when some predefined threshold is reached, maintenance is scheduled - \hspace{1sp}\citet{cigre445} elaborates on this in the context of power transformers. 
In the same context,~\citet{A_Mixed-Integer_Optimization_Model_for_Efficient_Power_Transformer_Maintenance_and_Operation} shows that these heuristics are inefficient in optimizing the profit of a single machine compared to a mixed-integer program (MIP). 

Research on optimizing maintenance at the precise moment when any delay in the maintenance action would result in the failure of a component is abundant. 
However, as shown in Example~\ref{ex:non_optimal_JIT_maintenance} from Appendix~\ref{sec:appendix}, this just-in-time maintenance strategy can lead to suboptimal or even arbitrarily bad solutions in our setting---undesirable when the underlying equipment is very costly.

Further, many of these methods resort to machine learning, using real-world data~\citep{Artificial_Intelligence_Application_for_Just_in_Time_Maintenance}.
As industries have only recently started collecting data systematically, approaches reliant on historical data are less reliable, making a mixed-integer nonlinear program (MINLP) more attractive in some settings.
Even in the presence of data, however, the reliability of MIP and its feasibility and optimality guarantees, as well as the knowledge of how far one is from the optimal solution, are also advantages over other artificial intelligence approaches.

\subsubsection{Production-maintenance scheduling} The problem discussed in this work belongs to the class of production-maintenance scheduling problems with parallel machines - see~\citet{Production_maintenance_and_resource_scheduling_A_review} for a comprehensive literature review. 
Many variants exist and have been studied for decades, albeit less than other more classical problems. 
The vast majority of the literature focuses on linear variants of this problem (see~\citet{Scheduling_jobs_with_time_resource_tradeoff_via_nonlinear_programming} for a resource-production exception), mostly because physical considerations tend to be ignored in more abstract models.
 Other authors avoid the nonlinearities by modeling machine failure with a random distribution~\citep{A_decomposition_method_by_interaction_prediction_for_the_optimization_of_maintenance_scheduling, Multi-component_Maintenance_Optimization:_A_Stochastic_Programming_Approach, A_Joint_Chance-Constrained_Stochastic_Programming_Approach_for_the_Integrated_Predictive_Maintenance_and_Operations_Scheduling_Problem_in_Power_Systems}, sometimes using Markov Chains to model a discrete set of conditions~\citep{Integrated_maintenance_planning_and_production_scheduling_with_Markovian_deteriorating_machine_conditions, markovian_model_liang_2018, Combined_Maintenance_Scheduling_and_Production_Optimization}.
  In these stochastic models, a few authors perform Bender's decomposition~\citep{Integrated_production_and_maintenance_scheduling_for_a_single_degrading_machine_with_deterioration-based_failures, Selected_topics_on_integrated_production-scheduling_and_maintenance-planning_problems}.

In his thesis,~\citet{Selected_topics_on_integrated_production-scheduling_and_maintenance-planning_problems} studies decomposition approaches for joint production-maintenance scheduling problems, using both Benders' and Dantzig-Wolfe decompositions.
This work considers the machines as single structures that are all identical.
 Machines are subject to random failures, and their degradation is implicitly modeled as the failure probability increases with the operational time.

In~\citet{A_Two-Stage_Framework_for_Power_Transformer_Asset_Maintenance_Management—Part_I:_Models_and_Formulations}, the authors heuristically decouple the maintenance and the production decisions to make the model more tractable.
 Not performing these decisions simultaneously allows the decision horizon for maintenance to be longer, but it loses reliability due to the assumptions on production. 

\subsection{Contributions}

In this article, we present an MINLP that models the production-maintenance scheduling problem for complex machinery.
To the best of our knowledge, our work is the first not to assume that the machines are inseparable units, but rather consider each as a structure composed of interacting components that can be independently maintained.
Furthermore, we allow nonlinearities in the degradation functions of the components, which results in MINLPs that can model real-world systems with greater accuracy at the cost of increased computational difficulty.
    
As solving this problem can be computationally expensive, we construct a Dantzig-Wolfe reformulation and solve the associated integer master problem with Branch-and-Price.
Branching is done on the master variables by adding hyperplanes to the master problem, complemented by a repair step to recover an integer solution in the original variable space.

  This article is organized as follows.
  In Section~\ref{sec:model}, a compact formulation is presented.
    Section~\ref{sec:DW} develops the Dantzig-Wolfe reformulation of this compact model, while Section~\ref{sec:branching_rule} presents the branching rule used to solve this reformulation to global optimality by a branch-and-price method.
    Implementation details and acceleration techniques appear in Section~\ref{sec:implementation}.
    Section~\ref{sec:tests} discusses the experimental setup and results.
    We conclude the paper in Section~\ref{sec:conclusion}, also indicating future research directions.

\section{Compact Formulation}
\label{sec:model}

In this production-maintenance scheduling problem, two types of decisions dictate the solution: how much to produce and when to maintain. 
The production, modeled by continuous variables, must satisfy the known demand at each period.
This impacts the auxiliary continuous variables modeling the condition of the components.
Since the condition of the components imposes a hard limit on production, this must be counterbalanced by component maintenance to guarantee satisfaction of the demand in subsequent periods.
Maintenance decisions are modeled by binary variables, and the minimization of their cost is the objective of this problem.
Maintenance actions are further restricted by the downtime constraints, which impose that no production occurs while maintenance is being performed.

Below, we present a general compact formulation of the production-maintenance scheduling problem described above.
Throughout this work, variables are represented by lowercase letters ($x,y,r$), parameters by uppercase letters in sans serif (\cmss{C}, \cmss{Q}, \cmss{M}, etc.), and sets by uppercase letters with calligraphic font ($\mathcal{N}, \mathcal{K}, \mathcal{I}$).
Subscripts $n,t$ will index machine $n$ and time period $t$, respectively, and superscript $k$ will refer to component $k$.
Tables~\ref{tab:abstract_parameter_description} and~\ref{tab:abstract_variable_description} present the parameters and variables of the compact formulation.

\begin{table}[!htb]
      \centering
          \resizebox{1\textwidth}{!}{
\begin{tabular}{ |p{1.5cm}||p{8cm}|p{3cm}|  }
     \hline
     Parameter & Parameter Meaning & Parameter Range\\
     \hline
    $\cmss{C}^k$   &  Cost of maintaining component $k$   & $\mathbb{R}^{+}$\\
    $\cmss{D}^k$  & Duration of component $k$'s maintenance & $\mathbb{N}$\\
    $\cmss{E}_t$   & Demand at period $t$ & $\mathbb{R}^{+}$\\
    $\mathcal{I}$  & Set of maintenance implications & $\mathcal{K} \times \mathcal{K}$\\
    $\mathcal{K}$   &  Index set of components to be maintained & $\mathcal{N}$\\
    \cmss{M} & Big constant to dominate constraint & $\mathbb{R}^{+}$\\
    $\mathcal{N}$ & Set of machines & $\mathbb{N}$\\
    $\cmss{Q}^k$   &  Maximum permissible production & $\mathbb{R}^{+}$\\
    $\cmss{R}^k$   &  Maximum condition of component $k$  & $\mathbb{R}^{+}$\\
    $\mathcal{T}$   &  Set of periods   & $\mathbb{N}$\\
    \hline
    \end{tabular}
    }
    \vspace{0.5em}
      \caption{Description of the parameters used in the compact formulation}
    \label{tab:abstract_parameter_description}
\end{table}


\begin{table}[!htb]
      \centering
          \resizebox{1\textwidth}{!}{

\begin{tabular}{ |p{1.25cm}||p{8cm}|p{2.5cm}|  }
     \hline
     Variable & Variable Meaning & Variable Type\\
     \hline
    $r^k_{n,t}$   &  Condition of component $k$ of machine $n$ at period $t$  & Continuous\\
    $x^k_{n,t}$   &  Maintenance of component $k$ of machine $n$ at period $t$   & Binary\\
    $y_{n,t}$   &  Production of machine $n$ at period $t$ & Continuous\\
    \hline
    \end{tabular}
}
        \vspace{0.5em}
        \caption{Description of the variables used in the compact formulation}
        \label{tab:abstract_variable_description}    
\end{table}


The production of machine $n$ at time $t$ is represented by the continuous decision variable $y_{n,t} \in [0, \min_{k \in \mathcal{K}^n} \cmss{Q}^k]$.
Here, $\cmss{Q}^k$ represents the production limit that component $k$ imposes when it is in perfect condition.
The requirement that the production over all machines must satisfy the demand is given by the following constraints: 

\renewcommand\theequation{1.\arabic{equation}}
\setcounter{equation}{0}
\begin{align}
  \sum_{n \in \mathcal{N}} y_{n,t} \geq \cmss{E}_t \quad \forall t \in \mathcal{T}\label{con:satisfy_demand}
\end{align}

The continuous variable $r_{n,t}^k \in [0,\cmss{R}^k_n]$ indicates the condition of component $k$ of machine $n$ at time period $t$.
Constraints~\eqref{con:production_limit} enforce a (possibly nonlinear) upper bound on the production of a machine, based on the condition of its components:
\begin{align}
  y_t \leq g_{n,k}(r_{n,t}^k) \quad \forall t \in \mathcal{T}\, \forall n \in \mathcal{N} \, \forall k \in \mathcal{K},{\label{con:production_limit}}
\end{align}
where $g$ is an increasing function.
The intent is that a machine with deteriorated components, and therefore a smaller value for $r_{n,t}^k$, will not be able to produce as much as a machine in better condition.

Binary decision variables $x^k_{n,t}$ model the maintenance decisions, indicating whether component $k$ of machine $n$ is maintained at time $t$.
Constraints~\eqref{con:maintenance_duration} impose a predetermined duration $\cmss{D}^k$ for maintenance:
 \begin{align}
  x^k_{n,i} \geq x^k_{n,t} - x^k_{n,t-1} \quad \forall n \in \mathcal{N} \, \forall k \in \mathcal{K}_n \, \forall t+1 \leq i \leq t + \cmss{D}^k{\label{con:maintenance_duration}}
 \end{align}

 Constraints~\eqref{con:maintenance_implications} represent maintenance action implications, where the replacement of some components forces the maintenance of others:
 \begin{align}
  x_{n,t}^k \leq x_{n,t}^{k^\prime} \quad \forall t \in \mathcal{T} \, \forall n \in \mathcal{N} \, \forall (k,k^\prime) \in \mathcal{I}_n {\label{con:maintenance_implications}},
 \end{align}
 where $\mathcal{I}_n$ is a given set of maintenance implications, meaning that if $(k, k^\prime) \in \mathcal{I}_n$, then the maintenance of component $k$ requires the maintenance of component $k^\prime$.

 Constraints~\eqref{con:maintenance_downtime} dictate that machines stop production while one of their components is being maintained:
 \begin{align}    
    y_{n,t} \leq (1-x^k_{n,t})\cdot \cmss{Q}^k \quad \forall t \in \mathcal{T} \, \forall n \in \mathcal{N} \, \forall k \in \mathcal{K}{\label{con:maintenance_downtime}}
 \end{align}

Finally, Constraints~\eqref{con:component_degradation} model the evolution of the condition of the components:
\begin{align}
  r_{n,t}^k \leq f_{n,k}(r_{n,t-1}^k, y_{n,t}; r_{n,t-1}^1, \dots, r_{n,t-1}^{|\mathcal{K}_n|}) + \cmss{M}\cdot x_{n,t}^k, \label{con:component_degradation}\\ 
  \forall t \in \mathcal{T} \forall n \in \mathcal{N} \, \forall k \in \mathcal{K}^n\nonumber{}
\end{align}

Function $f_{n,k}$ will model the (possibly nonlinear) degradation of component $k$ of machine $n$. 
It depends on the previous condition of the component, the current production, and may depend on the condition of other components of the same machine.
This work assumes that the components all start in optimal condition, so $r^k_{n,0} = \cmss{R}^k$, for all $n,k$.

The big-M term represents the maintenance action which, in conjunction with the binary variable $x_{n,t}^k$, will indicate the replacement of the component.
It does this by setting an upper bound on $r^k_{n,t}$ that is larger than $\cmss{R}^k_n$, thus allowing the component to be at its optimal state.
If function $f_{n,k}$ is concave, then the big-M is given by $f_{n,k}(0,\cmss{Q}^k; 0, \dots, 0) + \cmss{R}^k$.
It offsets the degradation by adding the biggest possible degradation and adds $\cmss{R}^k$, the component's optimal condition.

The objective of this problem is to minimize the maintenance cost (added over all components of all machines in all time periods). The complete compact model therefore is:

\begin{align}
\min_{x,y,r} \sum_{n \in \mathcal{N}}\sum_{t \in \mathcal{T}}\sum_{k \in \mathcal{K}^n} \cmss{C}^k \cdot x^k_{n,t}\\
\textrm{subject to: Constraints}~\ref{con:satisfy_demand}-~\ref{con:component_degradation}
  \label{model:abstract_formulation}
\end{align}

Given parameters $\mathcal{T}, \mathcal{N}$, and $\mathcal{K}_n, \mathcal{I}_n$ for all $n \in \mathcal{N}$, Table~\ref{tab:compact_formulation_size} presents the number of variables, constraints, and variable bounds of the model.

\begin{table}[H]
    \centering
    \begin{tabular}{|p{6em} | p{1em} p{15em}|}
        \hline
        Variables &  &$|\mathcal{T}|\sum_{n \in \mathcal{N}} (1+2|\mathcal{K}_n|)$\\[1em]
        Constraints & & $|\mathcal{T}|\sum_{n \in \mathcal{N}} (3|\mathcal{K}_n| + 1 + \frac{1}{D^k} + |\mathcal{I}_n|)$\\[1em]
        Bounds & & $6|\mathcal{T}|\sum_{n \in \mathcal{N}} (1 + |\mathcal{K}_n|)$\\
        \hline
    \end{tabular}
    \vspace{0.5em}
    \caption{Size of the model as a function of the input}
    \label{tab:compact_formulation_size}
\end{table}

The number of periods and the number of machines are the parameters with the most influence on the model size.
The machine's complexity, measured by the number of components and maintenance implications, further complicates things.

As this is a complex MINLP, valid constraints can help in solving it by removing parts of the search space known to be sub-optimal.
Since a maintenance action cannot start if it does not have enough time to finish, the following is a valid constraint:

\begin{align}
    x_{n,t}^k \leq x_{n,|\mathcal{T}|-\cmss{D}^k}^k, \forall k \, \forall t \geq |\mathcal{T}| - \cmss{D}^k
\end{align}

As will be seen in Section~\ref{sec:results}, this formulation is rather difficult to solve for bigger instances.
As such, we decided to explore a different approach to solve the problem.

\section{Extended Formulation}
\label{sec:DW}

The compact formulation's incidence matrix exhibits a block diagonal structure with a set of complicating constraints.
These are Constraints~\eqref{con:satisfy_demand}, requiring that the demand must be satisfied, while the blocks are composed of the constraints describing the functioning of the independent machines.
This structure is amenable to decomposition methods, given the near-independence of the block diagonal.
If the complicating constraints were to be removed, then the resulting model would be much easier to solve, as the remaining constraints can be partitioned by their index $n$ and solved independently.

Dantzig-Wolfe decomposition, originally introduced in~\citet{Decomposition_Principle_for_Linear_Programs}, can be used to reformulate problems with these attributes, so that they can be solved efficiently.
It begins by reformulating the problem into a pattern-based formulation, instead of the original assignment-based approach.
Assuming all feasible partial solutions satisfying each of the blocks are known, the formulation will decide which of these partial solutions should figure in the final solution. 
This formulation is called \textit{extended formulation}, alluding to its very high number of variables.
In contrast, the original formulation is called the \textit{compact formulation}.
Example~\ref{ex:formulations_example} below illustrates the difference between these two formulations with a simple abstract problem.
The example merely shows the difference in how the formulations view the problem, as the variables in the example are integer, to allow a more aesthetic presentation.

\begin{example}
  Consider the integer problem below, which has two independent blocks and a linking constraint.

\renewcommand\theequation{2.\arabic{equation}}
\setcounter{equation}{0}
  \begin{align}
    \min_{v,w} \quad      & v_1 + 2v_2 + w_1 + 3w_2 \\
    \textrm{subject to} \quad & v_1 + v_2 + w_1 + w_2 \geq 4 \label{con:ex_linking_con}\\
                        & 3v_1 + v_2 \leq 3 \label{con:ex_block1}\\
                        & 2w_1 + w_2 \leq 3 \label{con:ex_block2}\\
                        & v_1, v_2, w_1, w_2 \in \{0,1,2\}
  \end{align}

  The optimal solution is given by $v_1 = 0$, $v_2 = 2$, $w_1 = 1$, and $w_2 = 1$.
  Constraint~\ref{con:ex_linking_con} can be seen as a linking constraint, while Constraints~\ref{con:ex_block1} and~\ref{con:ex_block2} respectively
  represent the two independent blocks.
  Note how variables $v$ and $w$ only interact in Constraint~\ref{con:ex_linking_con}, and are otherwise independent.
  
  If we list all the ways of satisfying Constraint~\ref{con:ex_block1} and all the ways of satisfying Constraint~\ref{con:ex_block2}, then the optimal solution of the extended formulation corresponds to choosing the ways that satisfy Constraint~\ref{con:ex_linking_con} with minimum cost.
  Figure~\ref{fig:different_formulations} shows how the optimal solution would be represented in both formulations.

  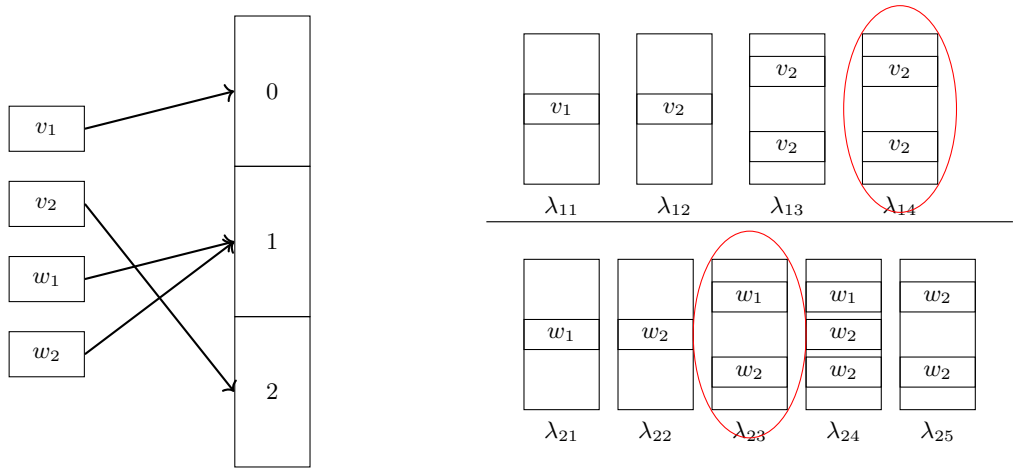
\begin{figure}[H]
    \centering
    \begin{minipage}{0.48\textwidth}
      \begin{tikzpicture}[
    box/.style={rectangle, draw, minimum width=1cm, minimum height=0.6cm, align=center},
    bin/.style={rectangle, draw, minimum width=1cm, minimum height=2cm, align=center},
    arr/.style={->, thick},
    font=\small
]

\node[box] (item1) at (0,4.5) {$v_1$};
\node[box] (item2) at (0,3.5) {$v_2$};
\node[box] (item3) at (0,2.5) {$w_1$};
\node[box] (item4) at (0,1.5) {$w_2$};

\node[bin] (bin1) at (3,5) {0};
\node[bin] (bin2) at (3,3) {1};
\node[bin] (bin3) at (3,1) {2};

\draw[arr] (item1.east) -- (bin1.west) node[midway, above] {};
\draw[arr] (item2.east) -- (bin3.west) node[midway, above] {};
\draw[arr] (item3.east) -- (bin2.west) node[midway, above] {};
\draw[arr] (item4.east) -- (bin2.west) node[midway, above] {};
\end{tikzpicture}
    \end{minipage}
    \hfill
    \begin{minipage}{0.48\textwidth}
      \begin{tikzpicture}[
    item/.style={rectangle, draw, minimum width=1cm, minimum height=0.4cm},
    bin/.style={rectangle, draw, minimum width=1cm, minimum height=2cm},
    font=\small
]

\node[bin] (bin11) at (0,0) {};
\node[bin] (bin12) at (1.5,0) {};
\node[bin] (bin13) at (3,0) {};
\node[bin] (bin14) at (4.5,0) {};

\node[bin] (bin21) at (0,-3) {};
\node[bin] (bin22) at (1.25,-3) {};
\node[bin] (bin23) at (2.5,-3) {};
\node[bin] (bin24) at (3.75,-3) {};
\node[bin] (bin25) at (5,-3) {};

\node[item] at (0, 0) {$v_1$};

\node[item] at (1.5, 0) {$v_2$};

\node[item] at (3, 0.5) {$v_2$};
\node[item] at (3,-0.5) {$v_2$};

\node[item] at (4.5, 0.5) {$v_2$};
\node[item] at (4.5,-0.5) {$v_2$};

\draw[line width=0.15mm] (-1, -1.5) -- (6, -1.5);

\node[item] at (0, -3) {$w_1$};

\node[item] at (1.25, -3) {$w_2$};

\node[item] at (2.5, -2.5) {$w_1$};
\node[item] at (2.5,-3.5) {$w_2$};

\node[item] at (3.75, -2.5) {$w_1$};
\node[item] at (3.75,-3) {$w_2$};
\node[item] at (3.75,-3.5) {$w_2$};

\node[item] at (5, -2.5) {$w_2$};
\node[item] at (5,-3.5) {$w_2$};

\node[below=1pt of bin11] {$\lambda_{11}$};
\node[below=1pt of bin12] {$\lambda_{12}$};
\node[below=1pt of bin13] {$\lambda_{13}$};
\node[below=1pt of bin14] {$\lambda_{14}$};

\node[below=1pt of bin21] {$\lambda_{21}$};
\node[below=1pt of bin22] {$\lambda_{22}$};
\node[below=1pt of bin23] {$\lambda_{23}$};
\node[below=1pt of bin24] {$\lambda_{24}$};
\node[below=1pt of bin25] {$\lambda_{25}$};

\node[ellipse, draw, red,minimum width=1.5cm, minimum height=2.75cm] at (4.5, 0) {};
\node[ellipse, draw, red,minimum width=1.5cm, minimum height=2.75cm] at (2.5, -3) {};

\end{tikzpicture}
    \end{minipage}
    \caption{Assignment-based and pattern-based approaches}
    \label{fig:different_formulations}
  \end{figure}

    In the optimal solution of the compact formulation, assignments are represented by arrows.
  In the optimal solution of the extended formulation, chosen patterns are represented by red circles.
  There is a clear equivalence between ways of assigning values to the variables of the compact formulation and ways of choosing patterns in the extended formulation.
  Notice also how each individual column in the pattern-based formulation satisfies the constraints associated with the respective subproblem.
For example, $\lambda_{13} = 1$ corresponds to $v_1=0, v_2=2$, which satisfies Constraint~\ref{con:ex_block1}.
Further, the columns chosen in the optimal solution together satisfy the linking constraints with the minimum cost.

The extended formulation is formally described below, where $\mathcal{Z}_i$ is the index set of the $\lambda$-variables of block $i$.
$v_{1j}$ and $v_{2j}$ (similarly for $w_{1j}$ and $w_{2j}$) represent the value of the original variable in $\lambda_{1j}$ (similarly $\lambda_{2j}$).
For example, $w_{14}=1, w_{24}=2$, which can be seen by looking at variable $\lambda_{24}$ in Figure~\ref{fig:different_formulations}.

\begin{align}
    \min_\lambda \quad & \sum_{j \in \mathcal{Z}_1} (v_{1j} + 2v_{2j}) \lambda_{1j} + \sum_{j \in \mathcal{Z}_2} (w_{1j} + 3w_{2j})\lambda_{2j} \\
    \textrm{subject to} \quad & \sum_{j \in \mathcal{Z}_i}\lambda_{ij} \leq 1, \forall i \in \{1,2\}\\
                              & \left(\sum_{j \in \mathcal{Z}_1} (v_{1j}+v_{2j})\lambda_{1_j}\right) + \left(\sum_{j \in \mathcal{Z}_2} (w_{1j}+w_{2j})\lambda_{2_j}\right) \geq 4\\
                & \lambda_{ij} \in \{0,1\}, \forall i \in \{1,2\}\,\forall j \in \mathcal{Z}_i
\end{align}

\xqed{$\triangle$}
\label{ex:formulations_example}
\end{example}

Only in toy problems is it possible to solve an extended formulation directly, since the number of feasible assignments tends to be huge.
Extended formulations are typically solved with column generation, a method that starts with a small subset of all the variables, and iteratively adds more, as long as they improve the previously obtained solution.
A key insight from LP theory is that there is always an optimal basic feasible solution, which has at most \textit{number-of-constraints} non-zero variables.
Since in an extended formulation the number of variables is so much larger than the number of constraints, this result indicates that relatively very few variables will need to be added to the model.
In the provided example, most of the columns did not need to be included.


There are then two types of problems.
The \textit{master problem}, which is responsible for coordinating the columns and guaranteeing that the complicating constraints are satisfied,
and the \textit{pricing problems}, which are responsible for generating the columns that must satisfy the respective block diagonal constraints.
These two problems are solved alternately, the master problem providing dual information which guides the pricing problems towards columns that will improve the current solution of the master problem. 

Column generation is a technique for solving linear problems, but the compact formulation in this work has both nonlinearities and binary variables.
While there are nonlinearities in the compact formulation, we delegate them to the pricing problems, thus not directly interfering with the column generation process. 
Examples of other works with nonlinear pricing problems can be found in~\citet{Exact_Methods_for_Recursive_Circle_Packing,Security_Constrained_Economic_Dispatch_Using_Nonlinear_Dantzig-Wolfe_Decomposition, Branch-and-price_for_a_class_of_nonconvex_mixed-integer_nonlinear_programs}.
In the presence of non-continuous variables, approaches employing column generation rely on a Branch-and-Bound tree as well.
At each node, the linear relaxation of the local problem is solved, and constraints are added to remove fractional solutions.
When the linear relaxation of every node in the Branch-and-Bound tree is solved with column generation, we have \textit{Branch-and-Price}.
For more details on column generation and Branch-and-Price see, e.g., ~\citet{Branch-and-Price}, or Uchoa et al.~\citet{Optimizing_with_Column_Generation:_Advanced_Branch-Cut-and-Price_Algorithms_Part_I}.

While nonlinearity will not impact our Branch-and-Price framework, the presence of continuous variables will.
Further explained in Section~\ref{sec:branching_rule}, they require special care when detecting integrality.
In~\citet{A_generic_view_of_Dantzig–Wolfe_decomposition_in_mixed_integer_programming}, the authors present a general framework for Dantzig-Wolfe decomposition in mixed-integer problems.
In their approach, they first split the columns into two, with one sub-column containing the values of the integer original variables and the other containing the values of the continuous original variables.
They later discretize the integer part, while enforcing convexity constraints on the continuous part.
Vanderbeck and Savelsbergh's approach was independently applied by~\citet{A_new_Dantzig-Wolfe_reformulation_and_branch-and-price_algorithm_for_the_capacitated_lot-sizing_problem_with_setup_times} for the Capacitated Lot Sizing Problem with Set Up Times.
Since our goal is to solve the problem described by the compact formulation, we allow the presence of relaxed integer variables in the optimal solution, as long as we can translate them into an integer solution in the original variable space.
For ease of exposition, these considerations will only be explicit when discussing the branching rule.

We will now start describing the application of these techniques to the problem we want to solve.

\subsection{Reformulating the production-maintenance scheduling problem}

To allow for an easier understanding, we present Tables~\ref{tab:extended_variable_description} and~\ref{tab:extended_parameter_description} with the variables and parameters for the extended formulation, respectively.

\begin{table}[!htb]
      \centering
          \resizebox{1\textwidth}{!}{

\begin{tabular}{ |p{1.25cm}||p{8cm}|p{2.5cm}|  }
     \hline
     Variable & Variable Meaning & Variable Type\\
     \hline
    $\delta_b$           & Variables satisfy all thresholds of branching constraint $b$ & \{0,1\} \\
    $\lambda_{z,i}$ & How many solutions $i$ are picked from subproblem $z$ & Integer \\
    $\theta$         & Dual variable associated to the convexity constraint & Continuous \\
    $\pi_t$         & Dual variable associated to time period $t$ & Continuous \\
    $\gamma_{b}$    & Dual variable associated to branching constraint b & Continuous \\
    \hline
    \end{tabular}
}
        \vspace{0.5em}
        \caption{Description of the variables used in the extended formulation}
        \label{tab:extended_variable_description}    
\end{table}

\begin{table}[!htb]
      \centering
          \resizebox{1\textwidth}{!}{
\begin{tabular}{ |p{1.5cm}||p{8cm}|p{3cm}|  }
     \hline
     Parameter & Parameter Meaning & Parameter Range\\
     \hline
    $\beta$            & Set of (variable, threshold) pairs for branching & $\mathcal{K}\times \mathbb{R}^+$ \\
    $\mathcal{B}$      & Set of branching decisions & - \\
    $\mathcal{L}_z$    & Index set of master variables for subproblem $z$ & $\mathbb{N}$ \\  
    $v_p$           & Threshold to split master variables on & $\mathbb{R}^+$ \\
    $\cmss{y}_{z,i,t}$ & Production of solution $i$ of subproblem $z$ at period $t$ & $\mathbb{R}^+$ \\
    $\mathcal{Z}$      & Set of aggregated subproblems & $\mathbb{N}$ \\
    $\cmss{Z}$         & Size of aggregated subproblems & $\mathbb{N}$\\
    \hline
    \end{tabular}
    }
    \vspace{0.5em}
      \caption{Description of the sets and parameters used in the extended formulation}
    \label{tab:extended_parameter_description}
\end{table}

In the following, for a given original variable $p$, the parameter $\cmss{p}_{i}$ denotes the value of $p$ in solution $i$.
So, for $y_t$ and $x^k_t$, we have $\cmss{y}_{i,t}$ and $\cmss{x}^k_{i,t}$. 
These use a different font to emphasize that they are parameters in the master problem.
$\mathcal{L}_n$ symbolizes the index set of master variables associated to machine $n$. 
In compact formulation terms, $i \in \mathcal{L}_n$ if the associated column satisfies all constraints for machine $n$, ignoring the demand constraints.
The subscript will be omitted when not relevant.
We now present the \textit{integer master problem} (IMP) below.

\renewcommand\theequation{3.\arabic{equation}}
\setcounter{equation}{0}
\begin{align}
\label{form:disaggregated_IMP}
\min_{\lambda}            &\sum_{n \in \mathcal{N}} \sum_{i \in \mathcal{L}_n} c^\intercal \lambda_{i,n}\\
\textrm{subject to}       &\sum_{n \in \mathcal{N}} \sum_{i \in \mathcal{L}_n} \cmss{y}_{i,t}\lambda_{i,n} \geq \cmss{E}_{t},  &  \forall t \in \mathcal{T} &\quad (\pi){\label{con:lambda_production_demand_dis}} \\
                          &\sum_{i \in \mathcal{L}_n} \lambda_{i,n} \leq 1, & \forall n \in \mathcal{N} & \quad (\theta){\label{con:lambda_convexity_dis}}\\
                          &\lambda_{i,n} \in \{0,1\}, & \forall n \in \mathcal{N} \, \forall i \in \mathcal{L}_n \label{con:integer_mvars_dis}
\end{align}

Constraints~\eqref{con:lambda_production_demand_dis} force the chosen solutions to satisfy the production demand, while Constraints~\eqref{con:lambda_convexity_dis} enforce that at most one solution per machine can be chosen.

The \textit{master problem} (MP) is obtained from the IMP by relaxing the integrality constraints. 
The master problem is the one that must be solved at every node of the Branch-and-Bound tree.

\subsection{Generating new columns}

Given that the master problem has a number of variables equal to the number of feasible partial assignments to each machine, it has an uncountably large number of variables, making it impossible to solve directly.
Thus, as is common with Dantzig-Wolfe decomposition, this extended formulation is solved with column generation.
The restricted master problem (RMP) sets all but a small subset of variables to $0$.
It starts with a restricted subset $\overline{\mathcal{L}} \subset \mathcal{L}$ of available columns. 
After solving the RMP, its optimal dual solution is used in a secondary problem, henceforth \textit{pricing problem}, tasked with finding columns that improve the optimal solution of the RMP.
The added column is the one in $\mathcal{L} \setminus \overline{\mathcal{L}}$ with minimum reduced cost.
The process is repeated until the column with the minimum reduced cost has a non-negative reduced cost, because adding it to the RMP would not improve the optimal solution.
It proves that none of the columns in $\mathcal{L} \setminus \overline{\mathcal{L}}$ would improve the optimal solution of the RMP, meaning it is also an optimal solution of the MP.

Variables $\pi$, lower bounded by $0$, represent the dual variables associated to the demand (Constraints~\ref{con:lambda_production_demand_dis}). 
Whenever precision is needed to indicate the dual variable referring to time period $t$, $\pi_t$ will be used.
Likewise, variables $\theta$, upper bounded by $0$, represent the dual variables associated to the convexity Constraints~\ref{con:lambda_convexity_dis}, and $\theta_n$ is the dual variable referring to machine $n$.
The reduced cost is then computed as $\cmss{C}^\intercal x - \pi^\intercal y - \theta_n$.
Since the production is a variable in this subproblem, $y$ uses the variable font. 

\renewcommand\theequation{4.\arabic{equation}}
\setcounter{equation}{0}
\begin{align}
\min_{x,y} \quad & \cmss{C}^\intercal x - \mathsf{\pi}^\intercal y - \theta_n\\
\text{subject to } \quad & \text{Constraints~\eqref{con:production_limit}-~\eqref{con:component_degradation} being satisfied}{\label{con:machine_operation}}
\end{align}

The dual variables $\pi$ incentivize solutions with high production at specific periods, and, informally, the dual variable $\theta$ indicates whether the solutions currently available for machine $n$ have a high impact on the objective of the solution of the RMP, making new ones valuable.

The value of the demand is a valid upper bound on the production of any specific machine, since any extra production would not be necessary, allowing us to add Constraint~\eqref{con:demand_upper_bound}.
\begin{align}
    y_t \leq \cmss{E}_t\, \forall t \in \mathcal{T}{\label{con:demand_upper_bound}}
\end{align}

Using the analysis done for the compact formulation in Section~\ref{sec:model}, the number of variables, constraints, and variable bounds of the pricing problem is presented in Table~\ref{tab:pricing_problem_size}.

\begin{table}[H]
    \centering
    \begin{tabular}{| p{6em} | p{1em}  p{12em}|}
        \hline
        Variables &  & $|\mathcal{T}| (1+2|\mathcal{K}_n|)$\\[1em]
        Constraints & & $|\mathcal{T}| (3|\mathcal{K}_n| + 1 + \frac{1}{D^k} + |\mathcal{I}_n|)$\\[1em]
        Bounds & & $6|\mathcal{T}| (1 + |\mathcal{K}_n|)$\\
        \hline
    \end{tabular}
    \vspace{0.5em}
    \caption{Size of pricing problem n as a function of the input}
    \label{tab:pricing_problem_size}
\end{table}

By being machine-specific, the pricing problem drops the sums over all machines of the compact formulation, as seen in Table~\ref{tab:compact_formulation_size}.
With many machines, this represents a large difference in size and difficulty.
The nonlinearity of the pricing problem will depend on the characteristics of $f$ and $g$ (Constraints~\ref{con:component_degradation} and \ref{con:production_limit}).

 The pricing problem can also be interpreted in the context of machine production. 
 The objective function is equivalent to the maximization of the profit of a machine (where the selling price of the product at time $t$ is given by the corresponding dual value). 
 The added valid inequality enforces a cap on production, after which the product cannot be sold, simulating markets where there is limited demand.
 An application of this reinterpreted problem into profit maximization of power transformers was introduced in~\citet{An_Optimization_Model_for_Power_Transformer_Maintenance}.

To prepare the exposition of the branching rule used in this work, we must first present the concept of identical subproblem aggregation.
Without aggregation, the branching rule would have required substantially less effort.
Instances of this production-maintenance scheduling problem may have identical machines, a common occurrence in real-world scenarios.
Whenever this is the case, the restricted master problem exhibits a lot of symmetry, requiring much additional effort.
To see this, suppose machines $n_1, n_2$ are identical. 
The pricing problem for machine $n_1$ is producing columns that also satisfy the non-linking constraints related to machine $n_2$, but are only available to machine $n_2$ after solving the pricing problem associated with it.

A common way to break these symmetries is to aggregate identical subproblems.
That is, given identical machines $n_1, n_2$, we will choose two columns out of one representative set instead of choosing one column each out of two identical sets. 
More formally, instead of choosing $\lambda_{n_1,i} \in \mathcal{L}_1$ and $\lambda_{n_2,j} \in \mathcal{L}_2$ with $\mathcal{L}_1 = \mathcal{L}_2$, we choose $\lambda_{n_1,i}, \lambda_{n_2,j} \in \mathcal{L}_1$. 

In the following, instead of variables $\lambda_{n,i}$ for each machine $n$, we will use variables $\lambda_{z,i}$ for each subgroup $z$ of identical machines, and the variables will indicate how many of these machines use column $i$. 
The parameter $\cmss{y}_{z,i,t}$ will denote the production in column $i$ from the subgroup $z$, at time $t$.
Further, the \textit{size of subgroup $z$}, the number of identical machines that were aggregated into subgroup $z$, is denoted by $\cmss{Z}_z$.
Thus, the variable $\lambda_{z,i}$ now has an upper bound equal to $\cmss{Z}_z$, rather than $1$, allowing us to choose the same column multiple times.
The aggregated integer master problem is described in~\ref{form:aggregated_IMP}. 

\begin{align}
\label{form:aggregated_IMP}
\min_{\lambda}            &\sum_{z \in \mathcal{Z}} \sum_{i \in \mathcal{L}_z} c^\intercal \lambda_{i,z}\\
\textrm{subject to}       &\sum_{z \in \mathcal{Z}} \sum_{i \in \mathcal{L}_n} \cmss{y}_{z,i,t}\lambda_{i,z} \geq \cmss{E}_{t},  &  \forall t \in \mathcal{T} &\quad (\pi){\label{con:lambda_production_demand}} \\
                          &\sum_{i \in \mathcal{L}_z} \lambda_{i,z} \leq \cmss{Z}_z, & \forall z \in \mathcal{Z} & \quad (\theta){\label{con:lambda_convexity}}\\
                          &\lambda_{i,z} \in \mathbb{Z}^+_0, & \forall z \in \mathcal{Z} \, \forall i \in \mathcal{L}_z \label{con:integer_mvars}
\end{align}

With this aggregation, we expect that the usefulness of the reformulation will increase with the number of identical machines, given its symmetry-breaking characteristics.
\section{An Exact Branch-and-Price Method}
\label{sec:branching_rule}

Up until now, the focus has been on the column generation algorithm, which solves LPs.
However, the compact formulation that we want to solve has integer variables.
This means that after the optimal solution of all pricing problems has a non-negative reduced cost, hence proving that the solution of the RMP is also optimal for the MP, branching is required.
As each node tends to be more computationally expensive in Branch-and-Price than traditional Branch-and-Bound, the branching rule gains increased importance.
In this section, we will use the notation $\lambda^\ast$ to indicate the value that variable $\lambda$ takes in an optimal solution, which we assume unique, for simplicity, unless specified otherwise.

Branching on single fractional master variables is not advisable in branch-and-price, as it leads to very unbalanced trees.
Given a fractional master variable $\lambda_i$, the added branching constraints would be $\lambda_i \geq \lceil \lambda_i^\ast \rceil$ and $\lambda_i \leq \lfloor \lambda_i^\ast \rfloor$.
The first constraint is extremely restrictive, forcing the column corresponding to $\lambda_i$ to be picked at least $\lambda_i \geq \lceil \lambda_i^\ast \rceil$ times.
On the other hand, the second constraint is very lax, as it limits one column out of uncountably many, and it would be possible to choose a similar column that would perform the same role in the RMP.
For that it would be enough, for example, to change one of the continuous variables by a small $\varepsilon$.
Furthermore, the branching constraints also need to be respected in the pricing problem to guarantee that this same pattern is not generated again, which would not be easy with this branching rule.

On the other hand, branching on pricing variables directly is not possible due to the aggregation of identical subproblems---two identical machines might pick different columns in the optimal solution.
While possible to disaggregate the subproblems as new original variable branching candidates appear, this would reintroduce symmetry to the problem, which we want to avoid.

For these reasons, we decided to branch on a set of bounds on the pricing variables as originally described in~\citet{On-Dantzig-Wolfe-decomposition-in-integer-programming-and-ways-to-perform-branching-in-a-branch-and-price-algorithm} and elaborated upon in~\citet[p. 470]{Branch-and-Price}.
Keeping with the notation of the latter, let $\mathcal{F}$ be the index set of the fractional master variables, $\beta$ the set of bounds on the pricing variables, $\mathcal{F}_\beta$ the index set of fractional master variables whose corresponding columns respect the bounds in $\beta$, and $v_\beta = \sum_{i \in \mathcal{F}_\beta} \lambda_i^\ast$.
To ease the exposition, we will assume that there is a single aggregated subproblem.
The work can be easily generalized to the case of different subproblems.

\subsection{Variable-bound branching}

After solving the RMP to optimality, we choose a pricing variable $p$, which can be any of the binary maintenance variables $x$.
Then, we pick a threshold $v_p$ and define $\beta = \{(p, v_p)\}$, and $\mathcal{F}_\beta = \{i \in \mathcal{F} \mid (p_i, v_{p_i}) \in \mathcal{R}\}$, where $\mathcal{R}$ is a relationship.
In this work, $(a,b) \in \mathcal{R}$ if $a \leq b$ or $a \geq b$, depending on the criteria explained below.
In reality, both options need to be considered to explore the solution space fully, but, for the sake of presentation, we omit this, as it would increase the density of the text for little benefit.
As a small detail, we use a priority queue to choose which subset $s$ to focus on next, and if its corresponding sum is integral, then the queue grows with two additional subsets, $s1, s2$, with $s_1=s \cup \{p_{i+1}\leq \lfloor v_{p_{i+1}}\rfloor\}$, and $s_2=s \cup \{p_{i+1}\geq \lceil v_{p_{i+1}}\rceil\}$.
For simplicity, we will assume in this explanation that all thresholds are an upper bound on the pricing variables, making $\mathcal{R}$ represent a less-than-or-equal inequality.
In other words, we will consider the fractional master variables whose coefficient for variable $p$ in its corresponding column is less-than-or-equal to $v_p$.
Using the notation introduced already, $\mathcal{F}_\beta = \{i \in \mathcal{F} \mid \cmss{p}_i \leq v\}$.
Then, the branching constraint at the down branch would be $\sum_{i \in \mathcal{L} \mid \cmss{p}_i \leq v_p} \lambda_i \leq \left\lfloor v_\beta \right\rfloor$.
This assumption will be dropped later when describing the pricing problem modifications.

For any integer solution, we must have $\mathcal{F} = \emptyset$ and $v_\beta = 0$, and this process works whenever $v_\beta$ is fractional.
However, it may happen that the set of fractional master variables is not empty, but $v_\beta \in \mathbb{Z}$, which does not allow us to branch.
For example, it would be the case if $\mathcal{F} = \{\lambda_1, \lambda_2\}$, $\cmss{p}_1 = \cmss{p}_2$, and $\lambda^\ast_1 = \lambda^\ast_2 = 0.5$.
In this case, we repeat the process by picking a different original variable $p^\prime$ and associated threshold $v_{p^\prime}$.
By increasing the restrictions in $\beta$, we remove elements from $\mathcal{F}_\beta$, which in the worst-case scenario is reduced to a single fractional variable, thus allowing us to branch.
Even in this scenario, the branching rule is not the same as the single-variable branching described at the beginning of this section.
This is because any variable added to the RMP that satisfies the restrictions in $\beta$ will later be added to the branching constraint.

When the patterns of the variables we want to branch on have no continuous variables, this process for finding a $\beta$ with a fractional $ v_\beta$ must terminate, given the finite number of original variables and master variables in the RMP.
To see this, simply choose $\lambda_i, \lambda_j \in \mathcal{F}_\beta$ and a pricing variable $\tilde{p}$ such that without loss of generality, $\tilde{\cmss{p}}_i < \tilde{\cmss{p}}_j$.
Such a pricing variable must exist, for otherwise there would be two identical columns in the RMP.
Then, the threshold can be $v_{\tilde{p}} = (\tilde{\cmss{p}}_i+\tilde{\cmss{p}}_j)/2$, and adding $(\tilde{p}, v_{\tilde{p}})$ to $\beta$ excludes at least one variable from $\mathcal{F}_\beta$, $\lambda_j$.
In fact, in the worst case, we need to add $log(|\mathcal{F}|)$ restrictions to $\beta$, as proven in~\citet[Prop. 3]{On-Dantzig-Wolfe-decomposition-in-integer-programming-and-ways-to-perform-branching-in-a-branch-and-price-algorithm}.

 So, finally, when $v_\beta \notin \mathbb{Z}$, the left branch is given by

 \renewcommand\theequation{5.\arabic{equation}}
\setcounter{equation}{0}
 \begin{equation}
     \sum_{i \in \mathcal{L}_\beta} \lambda^\ast_i \leq \left\lfloor v_\beta \right\rfloor
 \end{equation}

\begin{example}

  Figure~\ref{fig:variable_bound_branching} below shows the process of finding a suitable subset of fractional variables whose values in the optimal solution are fractional.
  To simplify the presentation, we abuse the notation when we describe the inequalities as $x_1 \geq 1$, since in reality it should be $\{i \mid x_{1,i} \geq 1\}$.
  After identifying such a subset, the figure also illustrates the addition of the branching constraint to the RMP and its impact on the MP.
  
  Suppose that the optimal solution of the RMP has five fractional variables, $\lambda_1, \lambda_2, \lambda_3, \lambda_4, \lambda_5$, with optimal solution values $0.5, 0.5, 0.3, 0.2, 0.5$, respectively, and that their corresponding columns are $\begin{bmatrix} 0 \\ 0 \\1 \end{bmatrix}, \begin{bmatrix} 0 \\ 1 \\0\end{bmatrix}, \begin{bmatrix} 1 \\ 0 \\0 \end{bmatrix}, \begin{bmatrix} 1 \\ 0 \\1 \end{bmatrix}, \text{ and } \begin{bmatrix} 1 \\ 1 \\0 \end{bmatrix}$.

Before starting to describe the application of the procedure to this example, we note that Figure~\ref{fig:variable_bound_branching} sacrifices some rigor in favor of visual clarity.
Specifically, the cuts are depicted as a simplified illustration of the separation of the master variables according to the specified bounds and should not be interpreted literally, since a faithful representation of the algorithm would not be possible in two dimensions.
The points in the top diagrams represent the variables in the space of fractional master variables, a subset of the RMP, which is the dashed rectangle in the bottom diagram.
The points in the RMP that are outside the inner rectangle are the master variables with an integer value in the optimal solution.
The larger outer rectangle in the bottom diagram represents the master problem, with the points (master variables) that have not yet been added to the RMP.
Representing these discrete sets as rectangles is also a deliberate simplification made to improve the readability of the figure.

  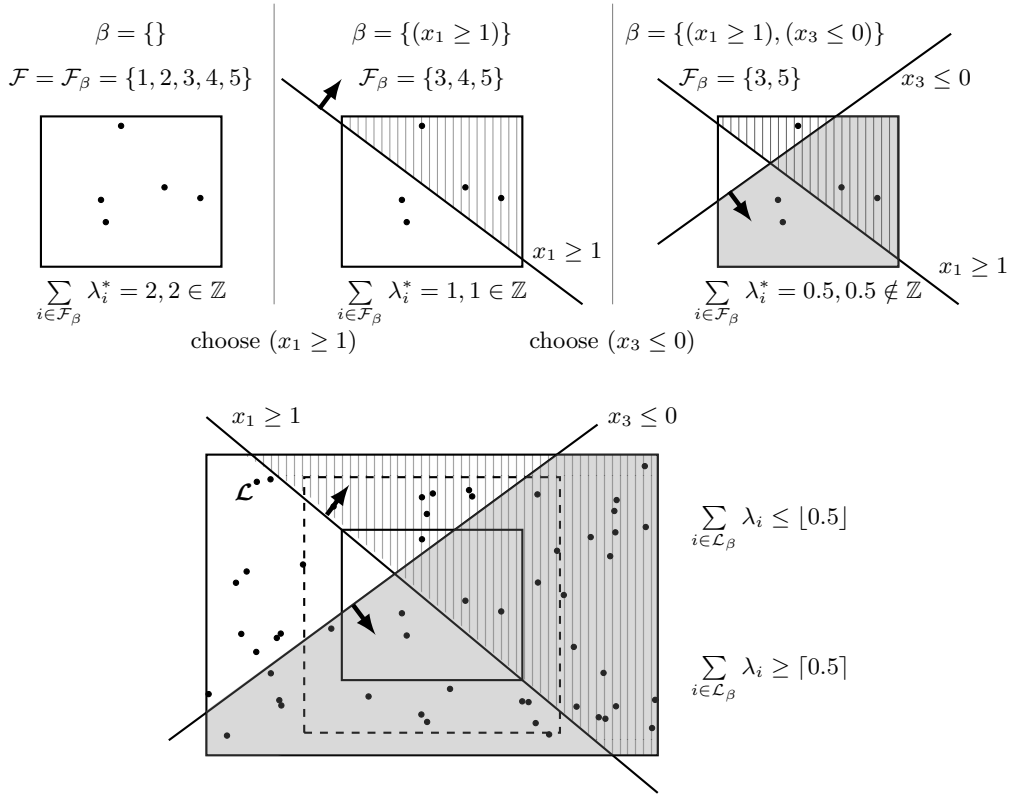
\begin{figure}[ht]
\centering

\begin{tikzpicture}[scale=1]
  \coordinate (center) at (0,0);

  \begin{scope}[shift={(-7,3)}]
    \node at (0,2.1) {$\beta = \{\}$};
    \node at (0,1.5) {$\mathcal{F} = \mathcal{F}_\beta = \{1,2,3,4,5\}$};
    \draw[thick] (-1.2,-1) rectangle (1.2,1);
    \pgfmathsetseed{37}
    \draw[scatter, only marks, mark=*, mark size=1pt]
    plot[domain=1:2, samples=2] ({-1.2 + 1*rnd}, {-1 + 1*rnd});
    \draw[scatter, only marks, mark=*, mark size=1pt]
    plot[domain=1:2, samples=3] ({-0.2+1.8*rnd}, {-0.2 + 1.4*rnd});
    \node at (0,-1.5) {$\sum\limits_{i \in \mathcal{F}_\beta} \lambda_i^* = 2, 2 \in \mathbb{Z}$};
  \end{scope}

  \draw[opacity=0.5] (-5.1,5.5) -- (-5.1,1.5);
\node at (-5.1,1) {choose $(x_1 \geq 1)$};

\begin{scope}[shift={(-3,3)}]
    \node at (0,2.1) {$\beta = \{(x_1\geq 1)\}$};
    \node at (0,1.5) {$\mathcal{F}_\beta = \{3,4,5\}$};
    \draw[thick] (-1.2,-1) rectangle (1.2,1);
    \draw[thick] (-2,1.5) -- (2,-1.5);
    \node at (1.8,-0.8) {$x_1 \geq 1$};
    \node at (0,-1.5) {$\sum\limits_{i \in \mathcal{F}_\beta} \lambda_i^* = 1, 1 \in \mathbb{Z}$};
    \pgfmathsetseed{37}
     \draw[scatter, only marks, mark=*, mark size=1pt]
    plot[domain=1:2, samples=2] ({-1.2 + 1*rnd}, {-1 + 1*rnd});
    \draw[scatter, only marks, mark=*, mark size=1pt]
    plot[domain=1:2, samples=3] ({-0.2+1.8*rnd}, {-0.2 + 1.4*rnd});

    \begin{scope}
      \clip (-1.2,-1) rectangle (1.2,1); 
      \fill[pattern=vertical lines, opacity=0.4] (-2,1.5) -- (2,-1.5) -- (2,2) -- (-2,2) -- cycle;
    \end{scope}
    \draw[line width=0.6mm, -latex] (-1.5,1.1) -- (-1.2,1.5); 
\end{scope}

  \draw[opacity=0.5] (-5.1+4.5,5.5) -- (-5.1+4.5,1.5);
\node at (-5.1+4.5,1) {choose $(x_3 \leq 0)$};

\begin{scope}[shift={(2,3)}]
    \node at (-0.7,2.1) {$\beta = \{(x_1\geq 1),(x_3\leq 0)\}$};
    \node at (-0.9,1.5) {$\mathcal{F}_\beta = \{3,5\}$};
    \draw[thick] (-1.2,-1) rectangle (1.2,1);
    \draw[thick] (-2,1.5) -- (2,-1.5); 
    \node at (2.2,-1) {$x_1 \geq 1$};
    \draw[thick] (-2,-0.7) -- (1.9,2.1); 
    \node at (1.7,1.5) {$x_3 \leq 0$};
    \node at (0,-1.5) {$\sum\limits_{i \in \mathcal{F}_\beta} \lambda_i^* = 0.5, 0.5 \notin \mathbb{Z}$};
        \pgfmathsetseed{37}
    \draw[scatter, only marks, mark=*, mark size=1pt]
    plot[domain=1:2, samples=2] ({-1.2 + 1*rnd}, {-1 + 1*rnd});
    \draw[scatter, only marks, mark=*, mark size=1pt]
    plot[domain=1:2, samples=3] ({-0.2+1.8*rnd}, {-0.2 + 1.4*rnd});
    \begin{scope}
      \clip (-1.2,-1) rectangle (1.2,1); 
      \fill[pattern=vertical lines, opacity=0.6] (-2,1.5) -- (2,-1.5) -- (2,2) -- (-2,2) -- cycle;
      \fill[gray, opacity=0.3](-2,-0.7) -- (1.9,2.1) -- (2,-2) -- (-2,-2) -- cycle;
    \end{scope}
    \draw[line width=0.6mm, -latex] (-1.05,0) -- (-0.75,-0.4);
\end{scope}

  \begin{scope}[shift={(-3,-2.5)}]
    \node at (4.5,1) {$\sum\limits_{i \in \mathcal{L}_\beta} \lambda_i \leq \lfloor 0.5 \rfloor$};
    \node at (4.5,-1) {$\sum\limits_{i \in \mathcal{L}_\beta} \lambda_i \geq \lceil 0.5 \rceil$};
    \draw[thick] (-3,-2) rectangle (3,2);
    \draw[dashed, thick] (-1.7,-1.7) rectangle (1.7,1.7);
    \draw[thick] (-1.2,-1) rectangle (1.2,1);
    \draw[thick] (-3,2.5) -- (3,-2.5);
    \node at (-2.2,2.5) {$x_1 \geq 1$};
    \draw[thick] (-3.5,-1.8) -- (2.2,2.4);
    \node at (2.8,2.5) {$x_3 \leq 0$};
    \node at (-2.5,1.5) {$\mathbfcal{L}$};
    \pgfmathsetseed{37}
    \draw[scatter, only marks, mark=*, mark size=1pt]
    plot[domain=1:2, samples=2] ({-1.2 + 1*rnd}, {-1 + 1*rnd});
    \draw[scatter, only marks, mark=*, mark size=1pt]
    plot[domain=1:2, samples=3] ({-0.2+1.8*rnd}, {-0.2 + 1.4*rnd});

    \draw[scatter, only marks, mark=*, mark size=1pt]
    plot[domain=1:2, samples=12] ({-3+1.8*rnd}, {-1.4 + 2*rnd});

    \draw[scatter, only marks, mark=*, mark size=1pt]
    plot[domain=1:2, samples=12] ({-3+6*rnd}, {1 + rnd});

    \draw[scatter, only marks, mark=*, mark size=1pt]
    plot[domain=1:2, samples=12] ({1.2+1.4*rnd}, {-2 + 4*rnd});

    \draw[scatter, only marks, mark=*, mark size=1pt]
    plot[domain=1:2, samples=12] ({-3 + 6*rnd}, {-1 - rnd});
    
    \begin{scope}
      \clip (-3,-2) rectangle (3,2);
        \fill[pattern=vertical lines, opacity=0.6] (-1.7,1.5) -- (3,-2.5) -- (3,2) -- (-2.5,2) -- cycle;
      \fill[gray, opacity=0.3] (-4.5,-2.5) -- (2.2,2.4) -- (20,-2) -- (-2,-2) -- cycle;
    \end{scope}
    \draw[line width=0.6mm, -latex] (-1.4,1.2) -- (-1.1,1.6);
    \draw[line width=0.6mm, -latex] (-1.05,0) -- (-0.75,-0.4);
  \end{scope}
\end{tikzpicture}

\caption{Visual representation of branching on original variable bounds. We use $x_1 \geq 1$ as shorthand for $\{i \in \mathcal{L} \mid x_{i,1} \geq 1\}$. 
The inner rectangle represents the discrete set of fractional variables in the optimal RMP solution, the dashed rectangle represents the discrete set of variables in the RMP, and the outer rectangle represents the discrete set of all the variables in the MP. The variables at the intersection of the two half-spaces are the ones in the new branching constraints.}
\label{fig:variable_bound_branching}
\end{figure}

    In the top left diagram, the sets are initialized.
    There are no bound changes yet, so $\beta$ is empty, and the set $\mathcal{F}$ of fractional master variables is the same as $\mathcal{F}_\beta$.
    The sum of all the fractional master variables satisfying the (empty) bound restrictions is $2$.
    Since this value is not fractional, the branching rule will not be able to create a branching constraint, and must keep refining the set $\mathcal{F}_\beta$.
    In the top middle diagram, the bound $x_1 \geq 1$ is chosen, and applying it to the previous set $\mathcal{F}_\beta$ gives us $\{3,4,5\}$, as the columns associated to $\lambda_1, \lambda_2$ do not satisfy this bound.
    Having reduced the number of fractional variables by nearly half, the resulting sum is $1$, still integer.
    This means that more bound changes are still required, and the introduction of $x_3 \leq 0$ in the top right removes variable $\lambda_4$.
    The resulting sum, $0.5$, allows us to branch.
    Since there is a finite set of possible bound changes, the process of finding a suitable set of bounds must always terminate.
    
    The set of master variables that satisfy the bound restrictions is then expanded with all the master variables in the RMP satisfying them, even those that were not fractional in the optimal solution, and those that were zero.
    Otherwise, the associated variables could have a negative reduced cost and be regenerated.

  Despite starting with only a limited view of the MP via the RMP, these restrictions are actually applied to the former.
  With this we mean that with the dynamic generation of columns, the branching constraints need to be updated throughout the iterations, whenever the generated column satisfies the bound restrictions.

  \xqed{$\triangle$}
\end{example}

Branching may be required multiple times throughout the optimization process, so we define $\mathcal{B}$ as the set of branching decisions that led to the current node.
Similarly, we define $\mathcal{F}_{\beta,b}$ as the index set of fractional master variables that satisfy the thresholds imposed by $\beta$ in branching decision $b \in \mathcal{B}$.
We will also drop the assumption on the thresholds, and define $\beta^\leq, \beta^\geq$ (similarly $\beta^\leq_b, \beta^\geq_b$) as the thresholds with a less-than-or-equal or greater-than-or-equal inequality, respectively.

\subsection{Incorporating dual variables in pricing}

As the branching decisions $b \in \mathcal{B}$ are added to the local RMPs, the corresponding duals, denoted by $\gamma_b$, will also influence the reduced cost of the remaining master variables. 
It is thus necessary to incorporate these duals into the pricing problem.
Since each branching decision $b$ restricts the master variables $\lambda_i$ in $\mathcal{F}_{\beta,b}$, the corresponding dual is non-zero only in these cases.
This justifies the introduction of binary variables $\delta_b$ for each branching decision $b$, which are $1$ if and only if all the thresholds in $b$ are respected, and consequently, whose dual variable should impact the objective.
Thus, we arrive at the pricing problem we must solve at every node.

  \begin{align}
    \min_{x,y} \quad &  \cmss{C}^\intercal x - \pi^\intercal y - \theta_z - \sum_{b \in \mathcal{B}_x} \gamma_b \delta_b\\
    \text{subject to} \quad & \left( \delta_b = 1 \Leftrightarrow (\forall p \in \beta^\leq_b ,\, p \leq v_p \land \forall p \in \beta^\geq_b ,\, p \geq v_p)\right), &\forall b \in \mathcal{B}  \label{con:branching_indicator}\\
                            & \delta_b \in \{0,1\},  &\forall b \in \mathcal{B} \\
                            & \text{Constraints~\ref{con:production_limit}-~\ref{con:component_degradation} being satisfied.} &
  \end{align}

Constraints~\ref{con:branching_indicator} are activated whenever all the thresholds relative to the branching constraint are satisfied.
We emphasize that these constraints neither forbid nor enforce any assignment of the pricing variables.
Rather, satisfying the thresholds of branching decision $\beta_b$ makes variable $\delta_b = 1$, resulting in the subtraction of the value of $\gamma_b$ from the reduced cost.
For a given branching decision $b$, the sign of $\gamma_b$ will either encourage or discourage the satisfaction of the thresholds of $\beta_b$.

These constraints can be linearized, but implementation becomes more complex in practice, since there needs to be a distinction between the up and down branches, as the indicator variables $\delta_b$ must be treated differently depending on the sign of $\gamma_b$.
The cases of negative and positive $\gamma_b$ also differ between the less-than-or-equal and the greater-than-or-equal thresholds.

If $\gamma_b$ is positive, then $\delta_b$ has an upward incentive, and we need to add constraints forcing it to $0$ if any of the thresholds are not satisfied. This can be translated to 

\begin{align}
    \delta_b \leq 1-p, \,&\forall p \in \beta_b^\leq \\
    \delta_b \leq p,   \,&\forall p \in \beta_b^\geq
\end{align}

On the other hand, if $\gamma_b$ is negative, setting $\delta_b$ to $1$ will increase the reduced cost.
So, we must force $\delta_b$ to $1$ if all thresholds are satisfied, which can be done with Constraint~\ref{con:force_delta_1}.

\begin{align}
    \delta_b \geq 1 - \sum_{p \in \beta_b^\leq} p +  \sum_{p \in \beta_b^\geq} p - \left (|\beta_p^\leq| + |\beta_p^\geq|\right )
    \label{con:force_delta_1}
\end{align}

\begin{example}
To show how the branching decision impacts the pricing problem, suppose the current node has depth $3$, and thus has two branching constraints.
As explained previously, these constraints have an associated set of thresholds to specify which master variables are affected.
In this example, for the first and second branching constraints, let the thresholds be $b_1 = \{(x_2 \leq 0)\}$ and $b_2=\{(x_1 \geq 1), (x_3 \leq 0)\}$, respectively.

Assume also that the current node is the down child of its parent, which in turn is the up child of the root node.
The branching dual is non-positive for the current branching constraint, and non-negative for the parent branching constraint.
See Figure~\ref{fig:gamma_tree_example} for a clarification.

\begin{figure}
\centering
\begin{tikzpicture}[every node/.style={circle,draw,minimum size=4mm},
    level distance=8mm, sibling distance=10mm]

  \node {\tiny Root}
    child {node {$\leq$}}
    child {node {$\geq$}
        child {node[fill=gray!40] {$\leq$}}
        child {node {$\geq$}}
      };

\end{tikzpicture}

\caption{The node in this example.}
\label{fig:gamma_tree_example}
\end{figure}
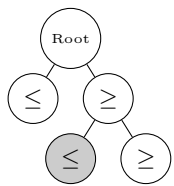

Figure~\ref{fig:gamma_activation} illustrates the behavior of the dual variables of the branching constraints in this example, as well as the incentives to satisfy the thresholds in the pricing problem.
From the sign of the dual branching variable, there is a pressure not to satisfy the thresholds of $b_1$, and another to satisfy the thresholds on $b_2$.
Depending on the remaining dual variables and constraints, the optimal solution of the pricing problem is more likely to satisfy $ \neg(x_2 \leq 0)$ and $(x_1 \geq 1, x_3 \leq 0)$.

\begin{figure}
\centering  

\begin{tikzpicture}[scale=0.8]
  \coordinate (center) at (0,0);

    \node at (-3,7) {$\mathcal{B} = \{b_1, b_2\}$};
    \node at (-3,6.5) {$b_1 = \{(x_2 \leq 0)\}\quad b_2 = \{(x_1 \geq 1), (x_3 \leq 0)\}$};

\begin{scope}[shift={(-7,3)}]
    \draw[thick] (0,0) ellipse (3 and 2);
    \draw[dashed, thick] (-1,2.5) -- (-1,-2.5);
    \node at (-2.2,2.5) {$x_2 \leq 0$};

    \node at (-2,0.1) {\large$\mathbf{\delta_{b_1} = 1}$};
    
    \begin{scope}
      \clip (0,0) ellipse (3 and 2);
        \fill[pattern=minus, pattern color=red, opacity=0.4] (-3,-2) -- (-3,2) -- (3,2) -- (3,-2) -- cycle;
        \fill[pattern=plus, pattern color=ForestGreen, opacity=0.4] (-1,-2) -- (3,-2.5) -- (3,2.5) -- (-1,2) -- cycle;
    \end{scope}
    \draw[line width=0.6mm, -latex] (-1.05,-0.8) -- (-1.6,-0.8);
  \end{scope}

\begin{scope}[shift={(2,3)}]
    \draw[thick] (0,0) ellipse (3 and 2);
    \draw[dashed, thick] (-3,2.5) -- (3,-2.5);
    \node at (-2.2,2.5) {$x_1 \geq 1$};
    \draw[dashed, thick] (-3.5,-1.8) -- (2.2,2.4);
    \node at (3,2.2) {$x_3 \leq 0$};

    \node at (1.5,0.1) {\large$\mathbf{\delta_{b_2} = 1}$};
    
    \begin{scope}
      \clip (0,0) ellipse (3 and 2);
          \fill[pattern=minus, pattern color=red, opacity=0.4] (-3,-2) -- (-3,2) -- (3,2) -- (3,-2) -- cycle;
        \fill[pattern=plus, pattern color=ForestGreen, opacity=0.4] (-0.4,0.4) -- (5,5) -- (5,-4) -- (-0.4,0.4) -- cycle;
        
    \end{scope}
    \draw[line width=0.6mm, -latex] (-1.4,1.2) -- (-1.1,1.6);
    \draw[line width=0.6mm, -latex] (-1.05,0) -- (-0.75,-0.4);
  \end{scope}

\end{tikzpicture}

\caption{Visual representation of the impact of
the branching decisions on pricing. The
ellipse represents the search space of the
pricing problem, and the patterns represent
the sign of the corresponding branching
variable $\gamma$. Only if the pricing variables
satisfy all the thresholds for a given
branching decision does the corresponding
dual affect the reduced cost. The thresholds
are dashed to indicate that these are
indicator constraints and that the thresholds
do not need to be respected.}
\label{fig:gamma_activation}

\end{figure}
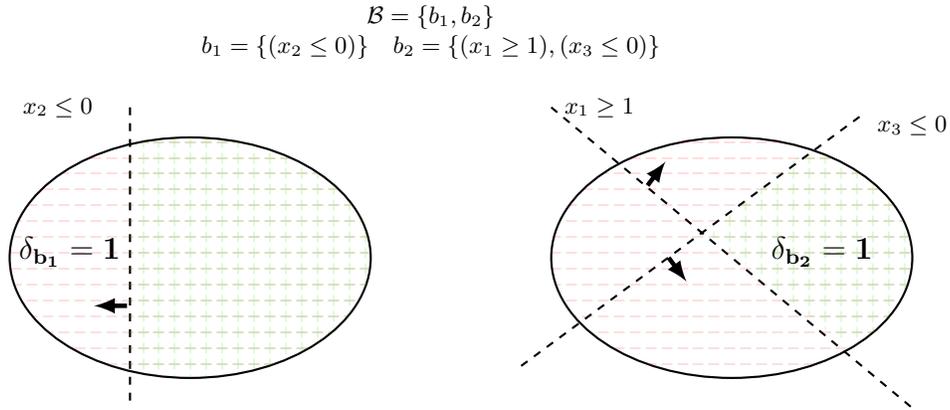

\xqed{$\triangle$}
\end{example}

From the works we have seen in the literature, variable-bound branching has always been achieved by finding hyperplanes that split the master variables into two fractional sets.
However, there is nothing forbidding the usage of nonlinear constraints for defining these sets.
That is, given original variables $p_1, \dots, p_n$ we can define $\mathcal{F}_\beta = \{i \in \mathcal{L} \mid h(\pp_{1,i},\dots, \pp_{n,i}) \leq 0\}$ for any function $h$, and the branching constraints in the pricing problem would look like $\delta_b = 1 \iff h_b(p_1, \dots, p_n) \leq 0, \forall b \in \mathcal{B}$.
Having nonlinear thresholds might more easily capture the fractional sets and allow one to rely on a solver's ability to find strong linearizations that might otherwise be too difficult to implement directly.
When the pricing problem is already an MINLP, this might achieve faster pricing more easily.
This approach will not be explored in this work. 

\subsection{Adapting the branching rule for the mixed-integer case}

At first sight, the branching rule we described does not seem to be complete, as illustrated by Example~\ref{ex:branching_rule_fractional} below.

\begin{example}
    Let $\mathcal{N} = \{1\}$, $\mathcal{K} = \{1\}$, and $\cmss{E}_1 = \cmss{E}_2 = 1$.
    Since there is only one machine and one component, we will omit the indices $n$ and $k$ in the notation.
    Further, let $f \equiv \cmss{R}$, and $g \equiv \cmss{Q} \geq \cmss{E}_1$.
    That is, there is no degradation, and the production limit is enough to satisfy the demand.
    Assume also that we have two master variables, $\lambda_1, \lambda_2$ with columns $\xx_{1,1}=\xx_{2,1}=0, \rr_{1,1} = \rr_{2,1} = \rr_{2,1} = \rr_{2,2} = \cmss{R}$, and $\yy_{1,1} = 2, \yy_{2,1} = 0, \yy_{1,2} = 0, \yy_{2,2} = 2$.
    In natural language, components are all at the optimal state, there is no maintenance, and columns $1$ and $2$ produce double the demand at time period $1$ and $2$, respectively, and produce nothing in period $2$ and $1$, also respectively.
    Then, the only optimal solution is given by $\lambda_1 = \lambda_2 = 0.5$.
    This solution is clearly optimal since it is feasible and its objective value is $0$.
    However, there is no hyperplane defined on original integer variables that can separate the fractional master variables into two disjoint sets, since $\xx_{1,1} = \xx_{2,1}$.
    Thus, this branching rule cannot find a hyperplane to separate the fractionalities.

    \xqed{$\triangle$}
    \label{ex:branching_rule_fractional}
\end{example}

We claim that whenever this is the case, there is an integer solution with the same objective value, which can be constructed by exploiting the convexity assumptions made in this work.
We apply such a procedure to the example above, but the explanation makes it clear that the approach is generalizable.

\begin{example}
    Under the conditions of Example~\ref{ex:branching_rule_fractional}, we add a new variable $\lambda_3$ such that, for every original variable $p$, we have $\cmss{p}_3 = \frac{\cmss{p}_1\lambda_1^\ast + \cmss{p}_2\lambda_2^\ast}{\lambda_1^\ast + \lambda_2^\ast}$.
    Thus, the new column is constructed by adding the coefficients of the columns associated with $\lambda_1,\lambda_2$, and dividing by their value in the original optimal solution.
    Since the compact formulation is convex, any convex combination of feasible solutions is also feasible.
    In particular, this statement remains true when fixing a machine index $n$, which corresponds to this convexification of the columns we just did.
    On the master problem side, since we start from a feasible solution, $\cmss{y}_{1,t}\lambda_1^\ast + \cmss{y}_{2,t}\lambda_2^\ast$ satisfies the demand for every period $t$.
    But then, so can $\lambda_3\frac{\cmss{y}_{1,t}\lambda_1^\ast + \cmss{y}_{2,t}\lambda_2^\ast}{\lambda_1^\ast + \lambda_2^\ast}$, if $\lambda_3 = \lambda^\ast_1 + \lambda^\ast_2$.
    With $\lambda_1=\lambda_2=0, \lambda_3=\lambda^\ast_1+\lambda^\ast_2$, the convexity constraints are also satisfied, by construction.
    Then, $\lambda_1 = 0, \lambda_2 = 0, \lambda_3 = 1$ is a feasible integer solution.
    Since the coefficients associated with the integer original variables of the master variables $\lambda_1, \lambda_2$ are the same (all zero), the objective value of this new solution is the same as that of the previous optimal solution.
    It is thus an integer optimal solution.
    
    \xqed{$\triangle$}
    \label{ex:branching_rule_repaired}
\end{example}

To get to this point, the sum of the master variables in the optimal solution must be an integer, so we can say that there is always an optimal solution where the variables we add take an integer value, and all other variables with the same maintenance pattern are $0$.

It is worth pointing out the need for normalizing by $\lambda^\ast_1+\lambda^\ast_2$ when constructing the new column, because it seems to be possible to avoid this, and then simply take $\lambda^\ast_1 = \lambda^\ast_2 = 0, \lambda^\ast_3 = 1$.
However, if the original fractional variables add up to an integer greater than one, the resulting column would not respect the non-linking constraints of the compact formulation---the convexity argument would no longer hold.
In other words, it is the difference between using a column $k$ times, and using ($k$ times a column) $1$ time.

This explanation was given under the assumption of a single aggregated subproblem.
If there are multiple subproblems, these arguments only hold if the $\lambda$-branching procedure cannot be applied in any of them.

Assuming we apply the repair step in the end, we can replace the original integrality constraints in Formulation~\ref{form:disaggregated_IMP} with alternative constraints that ensure that the repair step is possible, which we present next.

The theoretical portion of this paper will finish with the formal definition of the constraints we may use instead of the original integrality constraints~\ref{con:integer_mvars}.

Let us define an equivalence relation $\sim$ that groups master variable indices by the corresponding integer pattern.
In mathematical notation, for a given subproblem $z \in \mathcal{Z}$,

\begin{align*}
i \sim j \iff \cmss{x}_{i,t}^k = \cmss{x}^k_{j,t}, \forall k \in \mathcal{K}^z, \forall t \in \mathcal{T}
\end{align*}

This naturally defines equivalence classes $[j] := \{i \in \mathcal{L}_z \mid i \sim j\}$, and the set of all equivalence classes $\mathcal{L}_z/\sim$.
So, we replace the original constraints enforcing integrality with the following:

\begin{align}
    \sum_{i \in [j]} \lambda_i \in \mathbb{Z}, \forall z \in \mathcal{Z}\,\, \forall [j] \in \mathcal{L}_z/\sim \label{con:integrality_check}
\end{align}

In other words, all variables whose columns have the same coefficients for the maintenance variables must add up to an integer value.
As shown above, any solution satisfying these constraints can be converted to one satisfying the original integrality constraints, which can be done post-solve.

\begin{figure}[H]
    \scalebox{0.7}{
        \begin{tikzpicture}[node distance=2cm]
        
            \tikzstyle{startstop} = [ellipse, draw, text width=3.5cm, align=center]
            \tikzstyle{process} = [rectangle, draw, text width=4.5cm, align=center]
            \tikzstyle{decision} = [diamond, draw, text width=2cm, align=center]
            \tikzstyle{arrow} = [thick,->,>=stealth]

            \node (start) [startstop] {Node solved with\\ fractional solution $\lambda^\ast$ \\ $\mathcal{F}_\beta = \{i \in \mathcal{L}_z\}$};

            \node (decision1) [decision, below of=start, yshift=-2cm] {Constraint~\ref{con:integrality_check} satisfied?};

            \node (process1) [startstop, right of=decision1, xshift=4cm, text height=0.4cm, text width=2cm] {Solution optimal};

            \node (process2) [process, below of=decision1, yshift=-1cm] {Compute $v_\beta = \sum_{i \in \mathcal{F}_\beta} \lambda^\ast_i$};

            \node (decision2) [decision, below of=process2, yshift=-2cm] {$v_\beta$ \\fractional?};

            \node (process3) [startstop, right of=decision2, xshift=4cm] {Branch on \\$\sum_{i \in \mathcal{F}_\beta} \lambda_i \geq \lceil v_\beta \rceil$ or $\leq \lfloor v_\beta \rfloor$};

            \node (process4) [process, left of=decision2, xshift=-4cm, text width=4cm, text height=0.3cm]{Choose variable $p$ and threshold $v_p$};

            \node(process5) [process, above of=process4, yshift=-2cm, yshift=4cm, text width=4.5cm, text height=0.3cm]{Choose $\mathcal{R} \in \{\leq, \geq\}$};

            \node (process6) [process, above of=process5, yshift=2cm, text width=4cm, text height=0.3cm]{$\mathcal{F}_\beta = \mathcal{F}_\beta \cap \{i \mid (p_i, v_{p}) \in \mathcal{R}\}$};

            \draw [arrow] (start) -- (decision1);
            \draw [arrow] (decision1) -- node[anchor=south] {Yes} (process1);
            \draw [arrow] (decision1) -- node[anchor=west] {No} (process2);
            \draw [arrow] (process2)  -- (decision2);
            \draw [arrow] (decision2) -- node[anchor=south] {Yes} (process3);
            \draw [arrow] (decision2) -- node[anchor=south] {No} (process4);
            \draw [arrow] (process4.north) -- (process5.south);
            \draw [arrow] (process5.north) -- (process6.south);
            \path (process6.east) ++(1,0) coordinate (mid);
            \draw [arrow] (process6.east) -- (mid) |- (process2.west);
        \end{tikzpicture}
    }
    \caption{Branching on original variable bounds}
    \label{fig:branching_rule}
\end{figure}
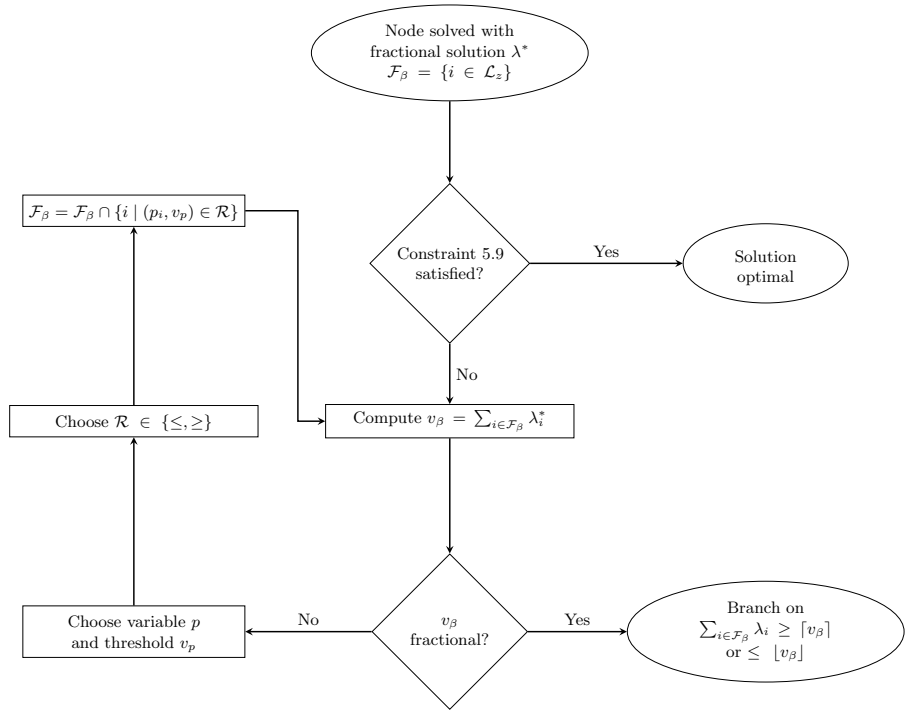

The branching rule is summarized in Figure~\ref{fig:branching_rule}. 

This concludes the theoretical part of the paper, and we will shift our focus to implementation details, acceleration techniques, and computational experiments.

\section{Implementation Details and Acceleration Techniques}
\label{sec:implementation}

For the implementation, we used SCIP~\citep{scip9.0}, as it is equipped with an easily customizable framework for Branch-and-Price. 
Nevertheless, there remained some problem-specific details that needed to be implemented, as described in this section.
We also highlight acceleration techniques for column generation and Branch-and-Price that are commonly found in the literature and were used in this work.

\subsubsection{Initial columns' generation}

In almost all instances, we use Farkas' Pricing~\cite[p.~43]{Generic_Branch-Cut-and-Price} for initializing the RMP and resolve infeasibilities during the Branch-and-Price process.
This technique uses Farkas' Lemma~\citep{Theorie_der_einfachen_Ungleichungen} on the RMP by computing patterns that break the lemma's infeasibility proof and adding them to the restricted master.
Implementation-wise, this corresponds to solving the pricing problem with the Farkas' duals and a zero cost coefficient.

We are not using Farkas' pricing in instances with many subproblems and few identical machines in each of them.
They are instead initialized with a modified run of the compact model.
This is done by running SCIP in the feasibility emphasis mode until a solution is found.

\subsubsection{Choice of subproblem}

As it may be expensive to solve the pricing problems exactly, we will not solve every subproblem in every iteration.
When a variable with a negative reduced cost is found, the RMP is solved again to get access to more updated duals.

We will use the results of a simple \textit{just-in-time} heuristic that fixes the production based on the dual values of the given iteration and tries to deduce the remaining variables, to decide the order in which the subproblems should be solved.
Subproblems with a more promising heuristic result will be solved first.

Besides this, pricing problems are initially solved with a gap limit and a solution limit.
Only when this limited run does not provide a variable with a negative reduced cost do we solve the problem to optimality.
Even with these restrictions, pricing took most of the solving time in our experiments, as can be seen in Section~\ref{sec:tests}.

\subsubsection{Choosing a branching candidate and associated threshold}

When deciding between a less-than-or-equal and a greater-than-or-equal inequality to define the fractional sets for branching, we choose the direction that maximizes the size of the respective set.
That is, if 

$$|\{i \in \mathcal{F}_\beta \mid \cmss{p}_i \leq \lfloor v_p \rfloor\}| > |\{i \in \mathcal{F}_\beta \mid \cmss{p}_i \geq \lfloor v_p \rfloor\}|,$$
then the less-than-or-equal inequality is chosen.
The greater-than-or-equal inequality is chosen if the second set is larger than the first.
The intuition is that more variables in the set, and consequently in the branching constraint, will lead to a more balanced tree.
In the case of a tie, we choose the inequality that allows us to keep the variables with a higher corresponding sum of squares, defaulting to less-than-or-equal.
That is, if 

$$\sum_{i \in \mathcal{F}_\beta \mid p_i \leq \lfloor v_p \rfloor} {\lambda_i^\ast}^2 \geq \sum_{i \in \mathcal{F}_\beta \mid p_i \geq \lceil v_p \rceil} {\lambda_i^\ast}^2,$$
then we opt for the less-than-or-equal inequality.
The reasoning is that these variables have a greater chance of being non-zero in the optimal solution. 

\subsubsection{Lower bounds}

Strong lower bounds are fundamental during the Branch-and-Bound process since they allow early pruning of suboptimal subtrees---no need to explore a subtree whose local lower bound is worse than the global primal bound.
The usual lower bound originating from the LP is not available since, at each pricing round, only the RMP is solved, not the MP.
It is thus necessary to employ alternative lower-bounding techniques.

The following, known as the Lagrangian dual bound, is a valid lower bound for the master problem~\cite[Prop 2.1]{Column_Generation}: 

\begin{equation}
c^\ast_{\text{RMP}} + \sum_{z \in \mathcal{Z}}\cmss{Z}_z\cdot \min(0,w_z),
\label{eq:lagrangian_bound}
\end{equation}

where $c^\ast_{\text{RMP}}$ is the optimal objective of the restricted master problem, and $w_z$ is a valid dual bound for pricing problem $z$.

The minimum reduced cost indicates the largest improvement in the objective per unit increase of the corresponding variable.
Thus, choosing to add the columns associated with the minimum reduced cost of each subproblem in the RMP yields an upper bound on the improvement of the current iteration. 
As not all the pricing problems are solved to optimality in every iteration, it is useful that~\ref{eq:lagrangian_bound} remains a lower bound to the master problem for any valid lower bounds on the pricing problems.

We can also use information from lower bounds of the pricing problem to stop the last pricing iteration earlier.
Given a primal bound $c_{\text{IMP}}$, we can stop pricing as soon as we find lower bounds $w_z$ for every pricing problem, such that~\ref{eq:lagrangian_bound} is not less than $c_{\text{IMP}}$.
While this trick does not reduce the number of pricing rounds, since the node would be pruned before the next pricing round, it lets one terminate the last pricing round sooner in every node that does not reach the branching stage.

In the case of a single subproblem $z$, the above can be transformed into a valid cut, as in every iteration, the reduced cost must be good enough to justify not pruning the node, for identical arguments.

\begin{equation}
    \cmss{C}^\intercal x - \pi^\intercal y - \theta_z - \sum_{b \in \mathcal{B}_z}\gamma_b\delta_b \leq \frac{c_{\text{IMP}}-c^\ast_{\text{RMP}}}{\cmss{Z}_z}
\end{equation}
In instances where the pricing problem is expensive, the valid cut can lead to a minor improvement.

If all subproblems are aggregated, because all machines are identical, then we also include the Farley bound, first discovered by~\citet{A_note_on_bounding_a_class_of_linear_programming_problems}.
For it, we need to define the \textit{minimum lambda price} $l(\pi)$:
\begin{equation}
l(\pi) := \min_{x \in \mathcal{L}} \left\{ \dfrac{x}{\pi^\intercal a_x} \bigg| \pi^\intercal_x > 0 \right\}.
\end{equation}

The minimum lambda price can also be used as an alternative to the minimum reduced cost for finding variables that would improve the RMP, which may be useful in some applications.
For more information, we refer to the ``Pivot rules and column management'' section of~\citet{Branch-and-Price}.
From this, the Farley bound can be derived:
\begin{equation}
l(\pi) \cdot \pi^\intercal b \leq c^\ast_{MP}, \quad \forall \pi \geq 0.
\end{equation}

Since the bound remains valid for every dual solution $\pi$, it can still be used even when using suboptimal dual values.
This is not the case with the Lagrangian bound, as the required lower bounds from the pricing problems would not be applicable.
Suboptimal dual values are often used in column generation, for example, when employing dual variable smoothing.
For information on dual variable smoothing, see~\citet{Weighted_Dantzig-Wolfe_decomposition_for_linear_mixed-integer_programming}.

The objective coefficients may also permit the tightening of these lower bounds.
For example, if all the variables in the master problem, even those not in the RMP, have objective coefficients that are multiples of $k$, then a bound $lb$ can be tightened to $k\cdot\left\lceil \frac{\text{lb}}{k} \right\rceil$, due to the integrality of the master variables of the IMP.
That is, if all coefficients are multiples of $1.5$ and we get a lower bound of $3.7$, then it can be tightened to $4.5$.

\subsubsection{Early branching} One of the common causes for a slow convergence of column generation is the tailing-off effect, where variables with near-zero reduced costs are generated.
Although many iterations can be spent generating such columns, this translates into very minor improvements to the relaxed master problem.
This effect can be lessened by a process called \textit{early branching}, where a node is not solved to optimality, but branching stops column generation when some condition is met---usually when the optimal reduced cost from the pricing problems gets close enough to zero.
The downside of early branching is that the optimal solution of the MP is not available, and thus cannot be used as a lower bound.
However, the lower bounds described above mitigate this issue, reducing the number of pricing rounds at the cost of potentially more branching decisions, depending on the tightness of the bounds.

When deciding to perform an early branching, there is a case that is of special interest.
When tightening the RMP solution yields a value not bigger than a known lower bound, we do not lose anything by performing early branching.
Mathematically, given a lower bound $lb$, if $\left\lceil c^\ast_{\text{RMP}} \right\rceil \leq \text{lb}$, solving the RMP to optimality will not result in a tighter integer bound.
In this case, early branching does not lose any relaxation strength.
This is based on the work in~\citet{Branch_and_price_for_the_length-constrained_cycle_partition_problem}.

\subsubsection{RMP heuristic} In comparison to Branch-and-Price, the application of column generation is straightforward, especially given the tools already available.
For this reason, and due to the strength of Dantzig-Wolfe relaxations, a common first attempt to determine the potential of an extended formulation is to solve the root node with column generation and then optimize the resulting RMP with the integrality restrictions.
This heuristic is called \textit{Price-and-Branch}, and is indeed a heuristic, since there is no guarantee that the optimal LP solution is integral.

This work is concerned with developing an exact method, but the Price-and-Branch heuristic remains quite useful for obtaining good primal bounds to quickly enable earlier pruning of subtrees.
We impose time limits on the solving of this integer problem, and a solution and objective limit.
As soon as a solution that improves the current primal bound is found, solving stops.
We use this heuristic after solving the root node, and again once enough new columns have been added since the last unsuccessful run of the integer restricted master problem.
\section{Computational Experiments}
\label{sec:tests}

 In this section, we will detail the experimental setup and metrics used to validate the models and assess their performance~\footnote{The code and experiment results can be found on GitHub~\citep{github-repo}.}.

 We generated a diverse set of instances randomly, parametrized by the following features:
 \begin{itemize}
     \item by the number of time periods (either $10$ or $20$),
     \item by the number of machines (either $20$ identical machines, or two groups of $10$ identical machines each), and 
     \item by their complexity (either ``high'' or ``low'').
\end{itemize}
 The complexity is a parameter of the machine that determines whether it will have a higher or lower number of components (between $1$ and $3$ and between $3$ and $7$, respectively) with more or less interaction between them (a $10\%$ chance of two specific components having a dependency, against $15\%$).
 The complexity of the machine also affects the amount of nonlinearity in the model, resulting from the increased inter-component interaction.
 For each of these combinations, we generated $50$ different instances.

 Functions $f$ and $g$ from the compact formulation and the pricing problem can randomly be linear, polynomial, or exponential, and the coefficients are sampled from a uniform distribution.
 Polynomial expressions are at most cubic, and
 exponential expressions use a fixed base and multiply the variable in the exponent by a constant.
 These constants are uniformly sampled between $0$ and $3$.
 Due to the nature of block-diagonal problems, exacerbated by the fact that not all components interact with each other, the instances are quite sparse.

 One of our first observations was that it is easy to create instances that favor each of the formulations, as having more or less identical machines favors the extended or the compact formulation, respectively. 
 For this reason, we opted for creating two sets of instances, with one of them having several identical machines (one group of $20$ machines), and another having two medium-sized groups of identical machines ($2$ sets of $10$ each).

The varying combinations of time periods, complexity, and machine groups yield $2\cdot 2\cdot 2\cdot 50 = 400$ different instances.
For each instance, we ran two experiments: the compact formulation and the extended formulation with exact pricing. They are referred to as ``Compact'' and ``DW'', respectively.

To obtain more meaningful results, we removed the instances that could be solved by both methods in under $5s$, and the ones for which no method could find an incumbent solution in the given time limit. 
This resulted in excluding $72$ and $107$ instances, respectively, and thus the results are shown for $221$ instances.
Of these, $118$ are feasible, and $103$ are infeasible.
We defend having this many infeasible instances because we believe that is likely to be how this work would be used most often in practice, with an industry expert tweaking the parameters (for example, the demand) to find a feasible and satisfactory model.
The results in the following section are also split between feasible and infeasible, since the methods performed very differently in these two subsets. 

In Table~\ref{tab:compact_dw} in Section~\ref{sec:results} below, the instances are split between ``easy'', ``medium'', and ``hard''.
Easy instances are those where one of the models could either prove optimality or infeasibility in under $10$s, medium instances where optimality or infeasibility was proven in the $300$s time limit by one of the models, and hard instances are those that hit the time limit and at least one of the models found a primal solution.
While the time and number of nodes shown are the arithmetic mean over the number of instances in the respective category, the gap is the arithmetic mean over the number of feasible instances in the category.
The gap is computed as $\frac{primal-dual}{\max(1,primal)}$, to deal with the cases where the dual bound is $0$.
In this table, the number of feasible instances is provided in parentheses to facilitate a better understanding of the methods' behavior.
For example, in the solved ``easy'' instances by the Compact formulation, we see $23(9)$.
This means that the Compact formulation solved $23$ easy instances, $9$ of which were feasible.

The experiments were run on a 6-Core Intel Core i7 with 64 GB running at 3.20GHz with a 300-second time limit.
The models were solved with the SCIP solver, version 9.2.1~\citep{scip9.0}.
The LP solver used was Soplex, version 7.1.3, and the NLP solver was IPOPT, version 3.14.17. 
For the implementation, we used SCIP's Python interface, PySCIPOpt~\citep{PySCIPOPT}, version 5.6.0, with Python version 3.13.3.

\subsection{Results and Analysis}
\label{sec:results}

Table~\ref{tab:compact_dw} below presents the results of the two approaches mentioned in this work across different difficulties, as explained in Section~\ref{sec:tests}.
We emphasize that the time and the number of nodes presented is the arithmetic mean over the number of instances, and the gap is the arithmetic mean over the number of feasible instances, defaulting to the primal bound when the dual is $0$.

\begin{table}[ht]
    \caption{Performance comparison of the different methods across different instance difficulties. In parentheses are the number of feasible instances.}
    \label{tab:compact_dw}
    \scriptsize
    
    \begin{tabular*}{\textwidth}{@{\extracolsep{\fill}} l c c c c c c c c c @{}}
    \toprule
    & \multicolumn{4}{c}{Compact} & \multicolumn{4}{c}{DW}\\
    \cmidrule(lr){3-6} \cmidrule(lr){7-10}
    Subset & instances & solved & time(s) & gap(\%) & nodes 
           & solved & time(s) & gap(\%) & nodes\\
    \midrule
    easy   & 57(9) & 23(9) & 189.1s & 5.6  & 475 & 57(9) & 6.1s & 0 & 5\\
    medium & 73(18) & 12(12) & 261.8s & 44.7 & 5238 & 62(7)  & 99.9s & 8.9 & 77\\
    hard   & 91(91) & -      & -      & 126.8 & 6426 & -      & -     & 11.3 & 47\\
    \bottomrule
    \end{tabular*}
\end{table}

Throughout the development process, a few results were observed consistently.
The decomposition approach, thanks to Farkas' Pricing, was very good at proving infeasibility, whereas the compact formulation struggled immensely with this.
The compact formulation, in turn, was better at proving optimality (21 instances vs 16), often due to the success of presolving.
In more challenging instances, the extended formulation was able to find good solutions much more quickly than the compact formulation, partly due to the RMP heuristic described in Section~\ref{sec:implementation}.

Focusing on the instances in the ``Hard'' subset, and comparing the results of the Compact and DW formulations, we see that the latter had a much smaller average gap ($126.8\%$ vs $11.3\%$), and again a much smaller number of nodes ($\sim 6.4$k vs $47$).
The number of nodes in particular showcases the strength of the branching rule used in this work.
However, while progress was consistently made on the compact formulation, some instances showed very little progress after the root node in the extended formulation, with the branching rule creating many nodes with the same local dual bounds.
Another observation is that the average number of nodes for the Dantzig-Wolfe decomposition was smaller in the ``Hard'' instances, when compared to the ``Medium'' instances.
This is because in these problems, the bottleneck was often the pricing problem, meaning that much time is spent solving the root node, and relatively few branchings are created.

\begin{figure}[H]
    \label{fig:DW_time_distribution}
    \centering
    \begin{tikzpicture}
    \def\printonlylargeenough#1#2{\unless\ifdim#2pt<#1pt\relax
    #2\printnumbertrue
    \else
    \printnumberfalse
    \fi}
    \newif\ifprintnumber
        \pie[
            text=legend,
            radius=1.8,
            color={blue!40, red!40, green!40, orange!40, purple!40},
            before number=\printonlylargeenough{2},
    after number=\ifprintnumber\%\fi
        ]{
            73/Exact pricing time,
            8/Python time,
            7/Integer RMP,
            11/Branching time,
            1/Others
        }
    \end{tikzpicture}

    \caption{Distribution of solving time across different components of the Dantzig-Wolfe decomposition, rounded to the nearest integer.}
    \label{fig:dw_dwpf_time}
\end{figure}
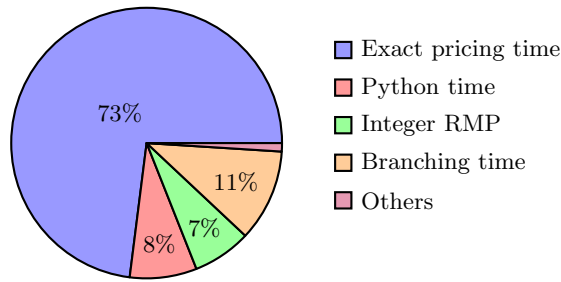

As expected, exact pricing was responsible for the vast majority of time spent, $72.82\%$, even with the limits imposed on it.
It makes a very strong case for the need for a heuristic to generate columns more quickly.
The branching rule is responsible for $11.35\%$ of the total solving time, suggesting that the algorithms used to choose variables and thresholds for the branching rule should be implemented in a more performant language.

Outside of SCIP, everything was implemented in Python.
As the data structures are quite complex and the codebase is big, the language's relatively slow performance led to a significant slowdown, making up $7.66\%$ of the total solving time.
The Price\&Branch heuristic, for obtaining master problem primal solutions, used $6.96\%$ of the total time.
While the time spent is not insignificant, we feel like the benefit of obtaining good primal solutions in these instances made up for it.
In any case, harsher solving limits can be put in place if desired.
The ``Others'' segment refers to the just-in-time heuristic for subproblem reordering (with $1.18\%$), and to the restricted master problem reoptimizations (with $0.03\%$).

\vspace{1em}

From this analysis, we conclude that the compact formulation yielded much worse results for difficult problems and for proving infeasibilities than the Dantzig-Wolfe decomposition.
On the other hand, it was slightly better at proving the optimality of easy instances than the decomposition approach.
In all, the extended formulation is the preferred method for solving this type of production-maintenance scheduling problems.

The decomposition would greatly benefit from a pricing heuristic, as exact pricing takes most of the time.
While that is a priority, rewriting the code in a faster language would also give some improvement.

\section{Conclusion and Future Work}
\label{sec:conclusion}

This article presents an MINLP for a general production-maintenance scheduling problem. 
Besides proposing a compact model, we also developed a Dantzig-Wolfe reformulation to take advantage of the block diagonal structure.

    The compact model can quickly solve smaller instances to optimality.
    However, it shows great difficulty in solving most instances, revealing a large gap in suboptimal instances and an inability to find incumbent solutions or prove infeasibility in most of them. 
    The decomposition approach was an improvement in almost every aspect, except for proving optimality, where the compact formulation proved slightly superior.
    Furthermore, the decomposition approach performed much better in harder instances, and also at detecting infeasibility, where the compact formulation yielded very poor results.
    Across all instance types, the extended formulation also used considerably fewer nodes.

  Regarding future work, it could be beneficial to explore the use of more intricate extended formulations with the goal of reducing the number of costly pricing rounds.
  Such formulations could be based on using separate maintenance and production patterns (patterns being an assignment to the corresponding variables) via two different easier pricing problems and combine them in the RMP. 
  These formulations would require more effort than the one presented in this paper, as it involves simultaneous column-and-row generation.

  On the compact model side, more than one type of maintenance action, which could differentiate between light and heavy maintenance, can make the model more interesting. 
  It would increase its difficulty, but it also makes it more realistic, as a complete replacement of a component is usually not the only option. 
  The model can also become more general. 
  For example, there might be non-critical components that may be maintained while still allowing for some production of the machine.
   Another example is the production upper bound in Constraints\eqref{con:production_limit}, which could depend not only on the condition of individual components but also on the joint condition of multiple components.
\appendix

\section{Appendix}
\label{sec:appendix}

This appendix provides an example showing that the \textit{just-in-time} maintenance strategy is only a heuristic.

\begin{example}
    Consider machines with two components, A and B, and suppose that the maintenance of component B implies the maintenance of component A.
Assume also that the production is heuristically planned independently of maintenance considerations.
Suppose that, by following this production plan, component A needs to be maintained in year $t_1$, while B can continue without maintenance until year $t_2 > t_1$.
Because maintenance of B implies maintenance of A, the \textit{just-in-time} heuristic here yields three maintenance actions in total.
Alternatively, in order to avoid the forced second maintenance of A, let us consider a solution where B and A are maintained both in year $t_1$.
There is nothing in our assumptions implying that this alternative solution is not feasible.
If it is feasible, then it would save us one maintenance of component A. 
For example, if the total number of periods is $t_2+1$, it is not unreasonable that both components could have been maintained at $t_1$ and still reached the end without failing.
This anticipation of a maintenance action would have resulted in a solution that is strictly better than the one we started with, proving that \textit{just-in-time} maintenance does not always produce the optimal solution.

\begin{figure}
    \centering
    \includegraphics[width=.5\textwidth]{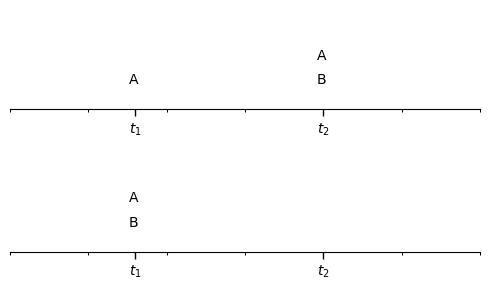}
    \caption{The anticipation of the B-A maintenance in the top solution leads to the bottom solution. 
    If feasible, the bottom solution is strictly better.}
    \label{fig:bad_greedy_solution}
\end{figure}

    Note also that if the bottom solution is feasible, then it is guaranteed to be an optimal solution for this fixed production vector. 
    By construction, the original solution tells us that component B requires one maintenance action at $t_2$ at the latest, and component A also requires one at $t_1$ at the latest.
    
    \xqed{$\triangle$}
\label{ex:non_optimal_JIT_maintenance}
\end{example}

The theoretical study of these local improvement steps is interesting but out of scope for the present paper.
Potential research directions are the detection of instances for which the local improvement step arrives at the optimal solution, and finding upper bounds on the number of steps to arrive at said optimal solution. 

\section*{Statements and Declarations}
        This work was funded by the Portuguese Foundation for Science and Technology - FCT with the scholarship 2023.00422.BD, and by Zuse Institute Berlin - ZIB.
        Part of the work for this article has been conducted in the Research Campus MODAL funded
        by the German Federal Ministry of Education and Research (BMBF) (fund numbers 05M14ZAM, 05M20ZBM)

\bibliographystyle{plainnat}
\bibliography{bibliography}

\end{document}